\documentclass[11pt]{article}

\makeatletter

\ifnum\@ptsize=0  \fi
\ifnum\@ptsize=2  \fi
\advance\textheight by 0.7cm
\advance\voffset by -0.3cm
\newcommand{\proof}{\noindent{\em Proof.}\, }
   \newcommand{\qed}{\hspace*{\fill}Q.E.D.\vskip12pt}

\newtheorem{theorem}{{\rm T\sc heorem}}[subsection]
    \newtheorem{lemma}[theorem]{{\rm L\sc emma}}
    \newtheorem{corollary}[theorem]{{\rm C\sc orollary}}
\newtheorem{proposition}[theorem]{{\rm P\sc roposition}}
\newtheorem{remark}[theorem]{{\rm R\sc emark}}

\def\rto{\cdots>}
\def\Pn#1{{\bf P}^{#1}}
\def\Pnd#1{{\check{\bf P}}^{#1}}
\def\Gr #1#2{{\bf Gr}(#1,#2)} 
\def\LG #1#2{{\bf LG}(#1,#2)} 
\def\OG #1#2{{\bf OG}(#1,#2)} 

\begin{document}
\title{
Geometry of the Lagrangian Grassmannian Sp(3)/U(3) with applications to
Brill-Noether Loci
}
\author{ Atanas Iliev\thanks{Partially supported by Grant MM-1106/2001
 of the Bulgarian Foundation for Scientific Research and by the 
 Norwegian Research Council}  \and  Kristian Ranestad   
}
\maketitle
\centerline{\it Dedicated to Andrei Nikolaevich Tyurin}
\vskip.5cm
\begin{abstract}
    The geometry of  $Sp(3)/U(3)$ as a subvariety of $\Gr 36$is explored to explain several examples given 
    by Mukai of non-abelian Brill-Noether loci, and to give some new examples. These examples identify 
    Brill-Noether loci of vector bundles on linear sections of the Lagrangian 
    Grassmannian $Sp(3)/U(3)$ 
    with orthogonal linear sections of the dual variety and vice versa. A main technical result of 
    independent interest is the fact that any nodal hyperplane section of the Lagrangian 
    Grassmannian projected from the node is a linear section of the 
    Grassmannian $\Gr 26$.

\end{abstract}
 
\section{Introduction}\label{sec:intro}

Mukai's linear section theorem for canonical curves of genus $8$ and $9$ says 
that every smooth canonical curve of genus $8$ with no $g^2_{7}$ is a complete intersection of the Grassmannian $\Gr 26\subset\Pn 
{14}$ with a 
linear subspace $\Pn 7$, and that every smooth canonical curve of genus 
$9$ with no $g^1_{5}$ is a  complete intersection of the Lagrangian Grassmannian
$\LG 36=Sp(3)/U(3)\subset \Pn {13}$ and a linear subspace $\Pn 8$ 
\cite{M4}.  Another beautiful theorem of Mukai is his interpretation 
of the general complete intersection of $\LG 36$ with a linear 
subspace $\Pn {10}$ as a non-abelian Brill-Noether locus of vector 
bundles on a plane quartic curve \cite{M3}.  This quartic curve can 
in a natural way be interpreted as the orthogonal plane section of 
the dual variety to $\LG 36$ in $\check\Pn {13}$.  In this paper we 
consider, very much in the spirit of \cite{M5}, general linear sections of the Lagrangian Grassmannian
$\LG 36\subset \Pn {13}$ of various dimensions, and show 
that the orthogonal linear section
of the dual variety $\check F$
of $\LG 36$  has an interpretation as a moduli space of
vector bundles  on the original linear section, and vice versa.  A 
similar study of linear sections of the $10$-dimensional spinor 
variety or orthogonal Grassmannian $\OG 5{10}\subset \Pn {15}$ is taken 
up by the first author and Markusewich in \cite{IM}.

The moduli spaces of stable vector bundles on curves is by now a classical 
subject dating back to 1960's and the fundamental work of Narasimhan, Seshadri 
and Tyurin (cf. \cite{NS}, \cite{T1}).  More recently the subvarieties of these moduli spaces 
representing bundles with many sections has attracted attention from  many 
authors (cf. \cite{B}, \cite{BF}, \cite{Go}, \cite{OPP}, \cite{Se}, 
\cite{M5}).  The 
corresponding theory for sheaves on $K3$ surfaces becomes 
particularly nice as explained in Mukais fundamental paper \cite{M84}.

The purpose of this paper is to present examples in this theory where 
these moduli spaces are all complete linear 
sections of either $\LG 36\subset\Pn {13}$ or its dual variety  of singular hyperplane 
sections in the dual space.   

A general tangent hyperplane section of $\LG 36$ is nodal, i.e. has a
unique tangency
point with a quadratic singularity.
The main technical result is (\ref{thm:hsection} and \ref{Zsurj})

\bigskip

\noindent
{{\rm T\sc heorem}}\label{theoremA}  {\it The projection of a nodal hyperplane
section of $\LG 36=Sp(3)/U(3)$ from the node
is a complete $5$-dimensional linear section of a Grassmannian variety $\Gr 26$.  This linear
section contains a $4$-dimensional quadric, and the general
$5$-dimensional linear section of $\Gr 26$ that contains a $4$-dimensional
quadric appears this way. }
\bigskip

We expect similar results to hold for the homogeneous varieties whose
general curve section are canonical curves of smaller genus.
In particular we expect the projection from the node of the general
nodal hyperplane section
of $\Gr 26$ to be a $7$-fold linear section of the spinor variety
$S_{10}$, and the projection from the node of the general nodal
$5$-fold linear section of $S_{10}$ to be the complete intersection of
a Grassmannian variety $\Gr 25$ and a quadric.  These lower genus 
cases will not be treated here.  

Given a smooth linear section $X$ of $\LG 36$ of dimension at most $4$,
each nodal hyperplane section that contains $X$ gives rise to an
embedding of $X$ into a $\Gr 26$.   In particular it gives rise to a
rank $2$ vector bundle on $X$
with $6$ global sections and determinant equal to the restriction of
the Pl\"ucker divisor on $\LG 36$.  Furthermore, this vector bundle is
stable, and when $X$ is at least $2$-dimensional, we
show that the set $\check F(X)$ of nodal hyperplane sections that contain $X$ form
a component of the corresponding moduli
space of stable rank 2 vector bundles on $X$ (cf. \ref{bnosurface}, 
\ref{bnothreefold}).  When $X$ is a curve, $\check F(X)$ form
a component of the corresponding Brill-Noether locus in the moduli
space of stable rank 2 vector bundles on $X$ (cf. \ref{curvecomp}).

In the opposite direction the general $3$-fold linear section $X$ of $\LG
36$ define a Brill-Noether locus of type II for the plane quartic curve
$\check F(X)$.  More precisely, associated to $X$ there is a $\Pn 1$-bundle
over $\check F(X)$ naturally embedded in $\Pn 5$ as a conic bundle of degree $16$.
It is isomorphic to a ${\bf P}({\cal F})$ where ${\cal F}$ is a rank $2$
vector bundle on $\check F(X)$ of degree $3$ such that no twist of ${\cal F}$ by
line bundles of negative degree have sections. Then $X$ is isomorphic to
the moduli space of rank two vector bundles ${\cal E}$ on $\check F(X)$ such that
rank$Hom_{{\bf C}}({\cal F,E})\geq 3$ and det${\cal E}-$det${\cal 
F}=K_{\check F(X)}$ (cf. 
\ref{fano threefold}).
\bigskip

The paper is organized as follows.  The first part is devoted
to the geometry of $\LG 36$.  It is the minimal orbit of an
irreducible representation of the symplectic group $Sp(3)$.  We
describe the four orbits of this group and analyse the singular
hyperplane sections corresponding to the corresponding $4$ orbits of
the group in the dual space.
$\LG 36\subset\Pn {13}$ parameterizes the $6$-fold of Lagrangian
planes in $\Pn 5$ with respect
to a given nondegenerate $2$-form.  A hyperplane section that is
singular at a point $p\in\LG 36$ define naturally a conic section in
the Lagrangian plane represented by $p$.  The correspondence between
singular hyperplane sections and conic sections in Lagrangian planes
manifests itself in various ways and is the crucial key to the main results of
this  paper.
In particular, the first part ends with a description of
conic bundles in the incidences between linear sections $X$, the
set $\check F(X)$ of singular hyperplane sections that contain $X$ and the
socalled vertex variety of points on conic sections in Lagrangian
planes corresponding to the singular hyperplane sections.

In the second part the conic section corresponding to a nodal
hyperplane section is again the crucial ingredient in the construction
of a rank $2$ vector bundle with $6$ global sections on the nodal
hyperplane section blown up
in the node.  The application to moduli spaces of vector bundles and 
Brill-Noether loci occupies the rest of
the second part and relies on some general results, to 
ensure that our examples form components of the corresponding 
moduli spaces.

\medskip\noindent
{\bf Notation.} We will denote by $\Gr kn$ the Grassmannian of rank 
$k$ subspaces of a $n$-dimensional vector space when $k<n$, and the 
Grassmanian of rank $n$ quotient spaces of a $k$-dimensional vector 
space when $n<k$.

\section{Geometry of the Lagrangian Grassmannian}\label{sec:geom}
\subsection{The groups $Sp(3,{\bf C})$ and
$Sp(3)$}\label{subsec:groups}
\paragraph{}\label{para:aa} Let $V = {\bf C}^{6}$ be a
$6$-dimensional complex vector space,
and let $\alpha : V \times V \rightarrow {\bf C}$,
${\alpha}: (v,v^{\prime}) \mapsto {\alpha}(v,v^{\prime})$ be a symplectic form
on $V$.
Thus ${\alpha}$ is bilinear, skew-symmetric and non-degenerate
(i.e. ${\alpha}(v \times V) = 0$ implies $v = 0$).

The canonical transformations of $V$
are those diffeomorphisms $f: V \rightarrow V$
which leave the form $\alpha$ invariant,
i.e. $f^* \ {\alpha} = \alpha$.
The canonical transformations of $V$ form a group
(with a multiplication -- the composition of maps),
and the symplectic group $Sp(3,{\bf C})$ is the
subgroup of these canonical transformations which
are complex-linear.

For the fixed $\alpha$, it is always possible to find
a base $\{e_1,...,e_{6}\}$ of $V$ in which
the Gramm matrix

\medskip

\centerline{$J = ({\alpha}(e_i,e_j)) =
\left( \begin{array}{cc}
0    & I_3 \\ -I_3 & 0
\end{array} \right)
$,}

\medskip

\noindent
where $I_3$ is the unit $3 \times 3$ matrix.
The complex symplectic group $Sp(3, {\bf C})$ has a natural embedding in
$GL(6, {\bf C})$ as the subgroup of all the complex rank $6$ matrices
which leave the matrix $J$ invariant, i.e.

\medskip

\centerline{$Sp(3,{\bf C}) = \{ Z \in GL(6, {\bf C}):{}^tZJZ = J \}$.}

\medskip

The group $Sp(3,{\bf C})$ is a non-compact Lie group
of real dimension $42$. Its Lie algebra $sp(3,{\bf C})$
consists of all the complex $6 \times 6$ matrices
of the form
\
$\left(
\begin{array}{cc}
A  & B \\
C & -{}^tA
\end{array} \right) $,
\
where $A,B,C$ are $3 \times 3$ matrices such that
${}^tB = B$ and ${}^tC = C$.
The group
$$
Sp(3,{\bf C}) = \{
\left(
\begin{array}{cc}
A & B \\
C &  D
\end{array} \right)
| B\cdot{}^tA = A\cdot{}^tB,\quad B\cdot{}^tC - A\cdot{}^tD =
-I_3,\quad D\cdot{}^tC = C\cdot{}^tD \}.
$$

\paragraph{}\label{para:ab}
Along with the non-compact group $Sp(3,{\bf C})$
there exists also the compact group $Sp(3)$,
which is defined below.

Let ${\bf R}$ be the field of real numbers,
and let
${\bf H} = {\bf R} + {\bf R}{\bf i} + {\bf R}{\bf j} + {\bf R}{\bf k}$
be the quaternionic algebra with products
${\bf i}^2 = {\bf j}^2 = {\bf k}^2 = -1$,
${\bf ij} = {\bf k} = -{\bf ji}$,
${\bf jk} = {\bf i} = -{\bf kj}$,
${\bf ki} = {\bf j} = - {\bf ik}$.

For the quaternion $h = a + b{\bf i} + c{\bf j} + d{\bf k}$,
the coordinate $a = Re(h)$ is the {\it real part} of $h$,
and
$\bar{h} = 2Re(h) - h = a - b{\bf i} - c{\bf j} - d{\bf k}$
is the {\it conjugate} quaternion to $h$.
The non-negative number $|h| = (a^2+b^2+c^2+d^2)^{1/2}$
is called the {\it norm} of $h = a+b{\bf i}+c{\bf j}+d{\bf k}$;
and it is easy to see that
$h\bar{h} = \bar{h}h = |h|^2$
for any $h \in {\bf H}$.

Let $e_1,e_2,e_3$ be a base of the $3$-dimensional quaternionic
vector space ${\bf H}^3 = {\bf H}e_1 + {\bf H}e_2 + {\bf H}e_3$.
In this base, any vector $u \in {\bf H}^3$ can be written
uniquely in the form

\medskip

\centerline{$u = {\Sigma}_{1 \le m \le 3} \ u_me_m
      = {\Sigma}_{1 \le m \le 3} \
                (a_m+b_m{\bf i}+c_m{\bf j}+d_m{\bf k})e_m$}

\centerline{= ${\Sigma}_{1 \le m \le 3} \
                           (a_m+b_m{\bf i})e_m + (c_m+d_m{\bf i}){\bf j}e_m$,}

\medskip

\noindent
where $u_m = a_m+b_m{\bf i}+c_m{\bf j}+d_m{\bf k}$,
$a_m,b_m,c_m,d_m \in {\bf R}$,
are the quaternionic coordinates of $u$, $m = 1,2,3$.

If $e_{3+m} = {\bf j}e_m, m = 1,2,3$,
then the above identity yields a natural identification of
${\bf H}^3$ and the complex $6$-space
$V = {\bf C}^{6} = {\Sigma}_{1 \le m \le 3} \ {\bf C}e_m + {\bf C}e_{3+m}$.

For any two vectors
$u = u_1e_1+u_2e_2+u_3e_3$ and $v = v_1e_1+v_2e_2+v_3e_3$
in ${\bf H}^3$, define their {\it quaternionic product}
$[u,v] = Re(u_1\bar{v}_1 + u_2\bar{v}_2 + u_3\bar{v}_3)$.
By definition

$$Sp(3) = \{ A \in GL(3,{\bf H})|
[Au,Av] = [u,v] \mbox{ for any } u,v \in {\bf H}^3 \}.$$

Clearly $Sp(3)$ is a subgroup of $GL(3,{\bf H})$;
and it is well-known that $Sp(3)$ is a compact connected
group of real dimension $21$.
In the above identification of ${\bf H}^3$ and ${\bf C}^6$,
any $A \in Sp(3)$ becomes unitary -- as an element of
$GL(6,{\bf C})$.
If $\alpha$ is the symplectic form on $V = {\bf C}^6$
with Gramm matrix $J$ as above, then
any $A \in Sp(n)$ becomes an element of $Sp(3,{\bf C})$.
Moreover $Sp(3) = U(6) \cap Sp(3, {\bf C})$,
where $U(6) \subset GL(6,{\bf C})$ is the unitary subgroup.

\medskip

The Lie algebra $sp(3)$ of $Sp(3) \subset GL(6,{\bf C})$
consists of all the complex $6 \times 6$ matrices
of the form
\
$\left(
\begin{array}{cc}
A  & B \\
-\bar{B} & \bar{A}
\end{array} \right) $,
\
where $A,B,C$ are $n \times n$ matrices such that
${}^tB = B$ and ${}^t{\bar{A}} = - A$.
The group
$$
Sp(3) = \{
\left(
\begin{array}{cc}
   A       &     B \\
-\bar{B} & \bar{A}
\end{array} \right)
| A\cdot{}^t\bar{A} + B\cdot{}^t\bar{B} = I_3,\quad B\cdot{}^tA =
A\cdot{}^tB \}.
$$

\subsection{The Lagrangian Grassmannian $\LG 36=Sp(3)/U(3)$}\label{subsec:LG}
\paragraph{}\label{par:ba}

The subspace $U \subset V$ is called {\it isotropic}
if ${\alpha}(U,U) = 0$.  The maximal dimension of an isotropic
subspace in $V$ is $3$, and in this case it is called {\it Lagrangian}.
Any isotropic subspace
is contained in some Lagrangian subspace.

By definition the complex Lagrangian Grassmannian
(or the complex Symplectic Grassmannian) $\LG 3V$
is the set of all the Lagrangian subspaces of $V = {\bf C}^6$.

The group $Sp(3, {\bf C})$ acts on the set of Lagrangian
     subspaces by $U\mapsto A\cdot U$, and it is easy to
     check that this action is transitive.

The Lagrangian Grassmannian $\LG 3V$
is a smooth complex $6$-fold
which admits a representation as a homogeneous space
$Sp(3,{\bf C})/St$ where $St$ is isomorphic to the stabiliser
group $St_U$ of any Lagrangian subspace $U \subset V$.
In our fixed base the subspace $U_0 = <e_1,e_2,e_3> \subset V$
is Lagrangian, and then
$St = St_{\scriptstyle U_0} \subset Sp(3,{\bf C})$
consists of all the matrices of the form
\
$\left( \begin{array}{cc}
A  & B \\
0 & {}^tA^{-1}
\end{array} \right) $,
\
where $A$ and $B$ are complex $3 \times 3$ matrices such that
$A.{}^tB = B.{}^tA$.
\paragraph{}\label{par:bb}
The Lagrangian Grassmannian has an alternative representation
as a quotient for the compact group
$Sp(3) = Sp(3,{\bf C}) \cap U(6) \subset Sp(3,{\bf C})$
by the subgroup $St_{\scriptstyle U_0} \cap U(6) \subset Sp(3)$.
It is easy to see that $St_{\scriptstyle U_0} \cap U(6)$ consists
of all the $6 \times 6$ matrices of the form
\
$\left( \begin{array}{cc}
A  & 0 \\
0 & {}^tA^{-1}
\end{array} \right)
$,
\
where $A \in U(3)$ is a unitary $3 \times 3$ matrix.
Now one can see that the Lagrangian Grassmannian
$\LG 3V$ is diffeomorphic to $Sp(3)/U(3)$
-- see \S 17 in \cite{TF};
from now on we shall use the former notation
or rather $\Sigma $ for its  Pl\"ucker embedding that we describe next.

\subsection{Representations and Pl\"ucker embedding}\label{subsec:rep}
\paragraph{}\label{par:ca}
  From here on we fix the form $\alpha$, the basis $\{e_1,...,e_{6}\}$
for $V$ and a dual basis $\{x_1,...,x_{6}\}$ for $V^{*}$.
With respect to this basis the matrix of the form $\alpha $ is
$$J = ({\alpha}(e_i,e_j)) =
\left( \begin{array}{cc}
0    & I_3 \\ -I_3 & 0
\end{array} \right)
,$$

\noindent
where $I_3$ is the unit $3 \times 3$ matrix. In the $\{x_{i}\wedge
x_{j}|1\leq i<j\leq 6\}$ for $\wedge^2V^{*}$,  the symplectic form
$\alpha$ has the  following expression

$$\alpha=x_{1}\wedge x_{4}+x_{2}\wedge x_{5}+x_{3}\wedge x_{6}.$$
It defines a correlation

$$L_\alpha :V\to V^{*}\quad v\mapsto\alpha(v,-).$$

This is natural up to sign and nonsingular since $\alpha$ is
nondegenerate.  It is called
a correlation since $v\in{\rm ker}L_\alpha(v)$ for every $v\in V$.
This correlation induces isomorphisms that we also denote by
$L_{\alpha}$:
$$L_{\alpha}:\wedge^kV \cong \wedge^kV^{*} $$
for $k=1\ldots 6$.
Consider the representations $\wedge^kV $  and $\wedge^kV^{*}$ of
$GL(6,{\bf C})$,
    and the induced representation of $Sp(3,{\bf C}) \subset GL(6,{\bf
    C})$.  The isomorphisms $\wedge^kV \cong \wedge^kV^{*} $ induced
    by $\alpha$ are clearly isomorphisms of these
     $Sp(3,{\bf C})$-representations.

    For the Pl\"ucker embedding of the Lagrangian Grassmannian we
    consider the $Sp(3,{\bf C})$ representation  $\wedge^3V$.  The form
    $\alpha$ defines a contraction which we denote by $\alpha$ itself:
    $\alpha:\wedge^3V\to V$.
    The representation then decomposes
    $$\wedge^3 V=V(14)\oplus V(6),$$
    where
    $$V(14)=\{w \in \wedge^3V| \alpha (w)=0\}$$
    and
    $$V(6)=\{w \in \wedge^3V| L_{\alpha}(w)\in \alpha\wedge V^{*}\}$$
   has rank $14$ and $6$ respectively (see e.g.
\cite{FH}, p.258).   Furthermore
$$V(14)^{*}:=L_{\alpha}(V(14))=\{\omega\in \wedge^3
V^{*}|\omega\wedge\alpha=0\}\subset \wedge^3
V^{*}$$

while $$V(6)^{*}=\alpha\wedge V^{*}.$$

   We now consider the decomposition of $V$ in two Lagrangian
subspaces $U_{0}=<e_{1}, e_{2}, e_{3}>$  and $U_{1}=<e_{4}, e_{5}, e_{6}>$.
We denote by $U_{0}^{\bot}$ the Lagrangian subspace
$L_{\alpha}(U_{0})=<x_{4},x_{5}, x_{6}>\subset V^{*}$, and likewise
$U_{1}^{\bot}=L_{\alpha}(U_{1})=<x_{1},x_{2}, x_{3}>$.
The decomposition
$$V^{*}=U_{1}^{\bot}\oplus U_{0}^{\bot},$$
induces a decomposition of $ \wedge^3{V^{*}}$:

$$ \wedge^3{V^{*}}=
\wedge^3U_{1}^{\bot}\oplus(\wedge^2U_{1}^{\bot}\otimes U_{0}^{\bot})\oplus
(U_{1}^{\bot}\otimes\wedge^2U_{0}^{\bot})\oplus \wedge^3U_{0}^{\bot}.$$

Consider the subspace
$$U_{1}^{\bot}\otimes\wedge^2U_{0}^{\bot}
\subset \wedge^3{V^{*}}.$$

The composition of the map $L_{\alpha}^{-1}: U_{0}^{\bot}\to U_{0}$, the
restriction $U_{1}^{\bot}\to {U_{0}}^{*}$ and the natural
isomorphism (up to scalar)
$\wedge^2U_{0}\to {U_{0}^{*}}$ define a natural map
$$U_{1}^{\bot}\otimes\wedge^2U_{0}^{\bot}\to {U_{0}^{*}}\otimes
\wedge^2U_{0}\to {U_{0}}^{*}\otimes {U_{0}}^{*}$$
Thus the component of a $3$-form in
$U_{1}^{\bot}\otimes\wedge^2U_{0}^{\bot}$ defines a bilinear form on
$U_{0}\subset V^{*}$.
Since the restriction $U_{1}^{\bot}\to {U_{0}}^{*}$ is an
isomorphism, the inverse map is a section  of the
restriction  $V^{*}\to {U_{0}}^{*}$.  Therefore the bilinear form is
independent of the choice of subspace $U_{1}$.

   We make this construction explicit in coordinates.
The exterior products
$e_{ijk} = e_i \wedge e_j \wedge e_k , 1 \le i < j < k \le 6$
    form a basis for $\wedge^3 V$.  Let $(x_{ijk}=x_i \wedge x_j \wedge
x_k)_{1 \le i,j,k \le 6}$
    be the dual basis.
   We interpret the basis for $\wedge^3V^{*}$ as the coordinate
   functions of an element
    $ w\in\wedge^3{V}$, and dually $\wedge^3{V}$ as
    coordinates on $\wedge^3V^{*}$.   Thus a $3$-form $\omega\in
    \wedge^3V^{*}$ has  coordinates:

\[
u^{*} = e_{123},
X^{*} = (x^{*}_{ab}) = \left(
\begin{array}{ccc}
e_{423} & e_{143} & e_{124} \\
e_{523} & e_{153} & e_{125} \\
e_{623} & e_{163} & e_{126}
\end{array}
\right) ,\]
\[
Y^{*} =
(y^{*}_{ab}) = \left(
\begin{array}{ccc}
e_{156} & e_{416} & e_{451} \\
e_{256} & e_{426} & e_{452} \\
e_{356} & e_{436} & e_{453}
\end{array}
\right) ,
z^{*} = e_{456}.
\]
The component of $\omega$ given by the matrix $X^{*}$ defines a bilinear form
on $U_{0}$.  If $(a_{1}, a_{2}, a_{3}), (b_{1}, b_{2},
b_{3})\in U_{0}$, then the form defined by
$X$ becomes:
$$(a_{1}, a_{2}, a_{3})\cdot \left(
\begin{array}{ccc}
e_{423} & e_{143} & e_{124} \\
e_{523} & e_{153} & e_{125} \\
e_{623} & e_{163} & e_{126}
\end{array}
\right)\cdot (b_{1}, b_{2},
b_{3})^t.$$ Similarly the component of $w$ given by the matrix
$Y$ defines a bilinear form
on $U_{1}$.
The subspace $V(14)^{*}\subset \wedge^3 V^{*}$ now has a simple
interpretation in these coordinates:
$$V(14)^{*} =
\{ \omega \in \wedge^3V^{*}: e_{i14} + e_{i25} + e_{i36} = 0, 1 \le i \le 6
\}.$$

In the decomposition
$\omega = [u^{*},X^{*},Y^{*},z^{*}]$ we get:

$$\omega \in V(14)^{*} \Leftrightarrow X^{*} = (x^{*}_{ab})\;\;
{\rm and}\;\;
Y^{*} = (y^{*}_{ab})\;\; {\rm are\; symmetric}\;\; 3 \times 3\;\;{\rm
matrices}$$

In particular, the bilinear form defined by the component $X^{*}$ is
the symmetric form
$$(x_{1}, x_{2}, x_{3})\cdot \left(
\begin{array}{ccc}
e_{423} & e_{143} & e_{124} \\
e_{523} & e_{153} & e_{125} \\
e_{623} & e_{163} & e_{126}
\end{array}
\right)\cdot (x_{1}, x_{2},
x_{3})^t=x^{*}_{11}x_{1}^2+\ldots + x^{*}_{33}x_{3}^2$$
   on $U_{0}$.   Thus we have described in coordinates a natural map
   $$q(U_{0}):\wedge^3U_{0}^{\bot}\oplus
U_{1}^{\bot}\otimes\wedge^2U_{0}^{\bot}\to Sym^2{U_{0}^{*}}$$

For $\omega\in \wedge^3U_{0}^{\bot}\oplus
U_{1}^{\bot}\otimes\wedge^2U_{0}^{\bot}$ we denote the associated
quadratic form
by $q_{\omega}(U_{0})$ or just $q_{\omega}$ if $U_{0}$ is understood
from the context, and by abuse we sometimes use the same notation for
the conic section in ${\bf P}(U_{0})$ that the quadratic form defines.
\medskip

\medskip

\paragraph{}\label{par:cb}

       The Pl\"ucker embedding $\Sigma := \LG 3V \subset
\Gr 3V \subset {\bf P}(\wedge^3V)$
of the Lagrangian Grassmannian $\LG 3V$
is the intersection of
$\Gr 3V$ with ${\bf P}(V(14))$, i.e.

$$\Sigma={\bf P}(V(14)) \cap \Gr 3V \subset {\bf P}(\wedge^3V).$$

\medskip

In coordinates the graph of a linear
map $X\in{\rm Hom}(U_{0},U_{1})$ is given by a matrix
$$
\pmatrix {I&X\cr}=\pmatrix {1&0&0&x_{11}&x_{12}&x_{13}\cr
0&1&0&x_{21}&x_{22}&x_{23}\cr 0&0&1&x_{31}&x_{32}&x_{33}\cr}.
$$
Thus
$$X\mapsto [1,X,\wedge^2X,\wedge^3X]$$
is an open immersion of ${\bf C}^9$ into $\Gr 3V$ in its Pl\"ucker
embedding defining an affine neighborhood of $[U_0]$.

This is the affine representation of coordinates
$[u,X,Y,z]$ dual to the coordinates $[u^{*},X^{*},Y^{*},z^{*}]$
above.  In particular, the subspace $U_X$ is Lagrangian if and only
if $X$ is
symmetric.  Thus the
intersection of the above coordinate chart of  $\Gr 3V$ with $\Sigma$ is
defined by $6$ symmetrizing linear equations of $X$ and $\wedge^2X$.
The closure of the affine chart around $[U_0]$ is clearly $\Sigma$ so
$\Sigma={\bf P}(V(14)) \cap \Gr 3V$ .

\medskip

\paragraph{}\label{par:dg}

Equations for $\Sigma$ in the affine coordinates
$$[1,X,\wedge^2X,\wedge^3X]=[1,x_{ij}, y_{ij},z]$$
\noindent
homogenized with first coordinate $u$, are given by standard
determinantal identities:

$$
\matrix{
uy_{11}-x_{22}x_{33}+x_{23}^2&uy_{12}+x_{12}x_{33}-x_{13}x_{23}\cr
uy_{22}-x_{11}x_{33}+x_{13}^2&uy_{13}-x_{12}x_{23}+x_{13}x_{22}\cr
uy_{33}-x_{11}x_{22}+x_{12}^2&uy_{23}+x_{11}x_{23}-x_{12}x_{13}\cr
uz-x_{11}y_{11}-x_{12}y_{12}-x_{13}y_{13}&zx_{11}-y_{22}y_{33}+y_{23}^2\cr
uz-x_{12}y_{12}-x_{22}y_{22}-x_{23}y_{23}&zx_{22}-y_{11}y_{33}+y_{13}^2\cr
uz-x_{13}y_{13}-x_{23}y_{23}-x_{33}y_{33}&zx_{33}-y_{11}y_{22}+y_{12}^2\cr
x_{11}y_{12}+x_{12}y_{22}+x_{13}y_{23}&zx_{12}+y_{12}y_{33}-y_{13}y_{23}\cr
x_{11}y_{13}+x_{12}y_{23}+x_{13}y_{33}&zx_{13}-y_{12}y_{23}+y_{13}y_{22}\cr
x_{12}y_{11}+x_{22}y_{12}+x_{23}y_{13}&zx_{23}+y_{11}y_{23}-y_{12}y_{13}\cr
x_{12}y_{13}+x_{22}y_{23}+x_{23}y_{33}&\cr
x_{13}y_{11}+x_{23}y_{12}+x_{33}y_{13}&\cr
x_{13}y_{12}+x_{23}y_{22}+x_{33}y_{23}&\cr
}
$$

\noindent
or, in short:
$$uY = {\wedge}^2 \ X, uzI = XY, zX = {\wedge}^2 \ Y,$$

\noindent
where $I = I_3$ is the unit $3 \times 3$ matrix. In other words,
the equations of $\Sigma \subset {\bf P}^{13}$ can be regarded
as projectivized {\it Cramer equations} for a symmetric
$3 \times 3$ matrix, its adjoint matrix and its determinant.

When $u=1$, the
equations define the coordinates $y_i$ and $z$ in terms of $x_i$, so
set-theoretically the coordinate chart $\{u=1\}$ is defined by these
equations.
So the
quadratic equations define a set-theoretic closure of this $U_0$.
An easy computation with MACAULAY shows that the resolution of
the ideal $I$ generated by these quadrics have betti numbers:

$$ \matrix{ 1 & - & - & - & - & -&-&- \cr
            - & 21 & 64 & 70 & - & -&-&- \cr
               - & - & - &-& 70& 64& 21& -  \cr
                 -& -& - & - & - & - & - & 1 \cr }$$

in the notation standard notation in the program (cf.\cite{MAC}).
Therefore the ideal $I$  is arithmetically Cohen-Macaulay and $Z(I)$
has no
embedded components.  Thus $Z(I)=\Sigma$.  In fact, one may check
that the resolution is symmetric, which
means that $I$ is arithmetically Gorenstein.

\subsection{Orbits, pivots and tangent spaces}\label{subsec:orbits}
\paragraph{}\label{par:da}
Recall from \cite{K} \S 9 that the action $\rho$ of the group $Sp(3, {\bf C})$
on $\Pn {13}={\bf P}(V(14))$ has precisely four orbits:

$$\Pn {13} \setminus F, F \setminus \Omega, \Omega \setminus \Sigma, \Sigma.$$

The dual action $\check{\rho}$ to is equivalent to $\rho$ induced by
$L_{\alpha}$,
so it has four corresponding orbits
$$\check{\bf P}^{13} \setminus \check{F}, \check{F} \setminus \check{\Omega},
\check{\Omega} \setminus \check{\Sigma}, \check{\Sigma}$$
in
$\check{\bf P}^{13} = {\bf P}({V}(14)^{*})$.
In this section we give a
   geometric characterisation of these orbits.

The smallest orbit, and the only closed one, is the Lagrangian
Grassmannian itself.  The closure of the orbit $F\setminus\Omega$ is the
union of the projective tangent spaces to $\Sigma$. They form a
hypersurface $F \subset {\bf P}^{13}$.  Similarly,
the dual variety to $\Sigma$, i.e.
the set of hyperplanes $H \subset {\bf P}^{13}$
containing some projective tangent space to $\Sigma$ form
a hypersurface $\check{F} \subset \check{\bf P}^{13}$ isomorphic to $F$.
By \cite{SK}, p. 108 the equation defining $F$ is
\begin{equation}\label{F}
f(w) = (uz - tr \ XY)^2 + 4u \ det \  Y + 4z \ det \ X
- 4{\Sigma}_{ij} \ det(X_{ij}).det(Y_{ij}),
\end{equation}
\noindent
where $X_{ij}$ and $Y_{ij}$ are the complimentary matrices
to the elements $x_{ij}$ and $y_{ij}$
(see also \cite{SK} p. 83).  $\Omega$ is the singular locus of $F$,
defined by vanishing all the partials of $f$. A simple computation
in MACAULAY \cite{MAC} shows that the ideal of $\Omega$ has a resolution with
Betti numbers
$$ \matrix{ 1 & - & - & -&- \cr
            - & - & - & - & -\cr
               - & 14 & 21 & - & -  \cr
                 -& -& - & 14 & 6  \cr }.$$
Furthermore $dim \ \Omega = 9$
and $deg \ \Omega = 21$ and sectional genus $22$.
We shall see this in a different description
of $\Omega$ later.

These orbits of $\rho$ in $\Pn {13} = {\bf P}(V(14))$
are restrictions to ${\bf P}^{13} = {\bf P}(V(14))$
of the orbits of $SL(V, {\bf C})$
in ${\bf P}^{19} = {\bf P}({\wedge}^3 \ V) =
{\bf P}(V(14) \oplus V(6))$.
In \cite{D}, Donagi proves the following:

\begin{theorem}\label{thm:segre} (Theorem of Segre for $\Gr 36$)
The natural representation of $GL(6,{\bf C})$
on ${\bf P}^{19}={\bf P}({\wedge}^3 \ V)$ has four orbits:

\medskip

\centerline{
${\bf P}^{19} = G \cup (W  \setminus G) \cup (D \setminus W) \cup 
({\bf P}^{19} \setminus D)$,
\ where:
}

\medskip

{\bf (i)} \
$G = \Gr 36$ is the Grassmannian $9$-fold of $3$-vectors
in $V = {\bf C}^6$, or

\smallskip

\centerline{$G = \{ \mbox{ 3-spaces } U: U \subset V \}$;}

\smallskip

{\bf (ii)} \
$D \subset {\bf P}^{19}$ is a quartic hypersurface
which is the union of all the tangent lines to $G$,
i.e.  $D = {\cup}_{u \in G}{\bf P}^9_u$,
where ${\bf P}^9_u$ is the tangent projective $9$-space to $G$
at $u \in G$;

\smallskip

{\bf (iii)} \
    $x\in {\bf P}^{19} \setminus W \Leftrightarrow$
    there exists a unique secant or tangent line $L_x$ to $G$ through $x$,
    and then $x \in D \setminus W$ \ {\it iff} \ $L_x$ is tangent
    to $G$;

\smallskip

{\bf (iv)} \
$W$ is the singular locus of $D$,
and $x\in W \Leftrightarrow$ there exist
infinitely many secant lines to $G$ passing through $x$.

\end{theorem}
\proof cf. \cite{D}, Ch. 3.
\qed

By restriction to ${\bf P}(V(14))\subset {\bf P}({\wedge}^3 \ V)$ we
get the analogous result for the orbits of $\rho$:

\begin{corollary}\label{cor:segreLG} (Theorem of Segre for $Sp(3)/U(3)$)

The representation ${\rho}$ of $Sp(3,{\bf C})$ on ${\bf P}^{13}={\bf P}(V(14))$
   has four orbits:

\medskip

\centerline{
${\bf P}^{14} =
\Sigma \cup ({\Omega}\setminus{\Sigma}) \cup (F\setminus{\Omega}) \cup
({\bf P}^{13}\setminus F)$,
\ where:}

\medskip

{\bf (i)} \
$\Sigma$ is the 6-fold Lagrangian Grassmannian

\smallskip

\centerline{ $Sp(3)/U(3) =
\{ \mbox{ Lagrangian 3-spaces } U: U  \subset V = {\bf C}^6 \}$
}

\smallskip

{\bf (ii)} \
$F \subset {\bf P}^{13}$ \ is a quartic hypersurface
which is the union of all the tangent lines to $\Sigma$,
i.e. $F = {\cup}_{u \in {\Sigma}}{\bf P}^6_u$,
where ${\bf P}^6_u = {\bf T}_u  \Sigma$
is the tangent projective $6$-space at $u \in \Sigma$;

\smallskip

{\bf (iii)} \
$w \in {\Pn {13}} \setminus \Omega$ \ $\Leftrightarrow$ \
there exists a unique secant or tangent line $L_w$ to $\Sigma$
through $w$,
and then $w \in F \setminus \Omega$ \ {\it iff} \ $L_w$ is tangent
to $\Sigma$;

\smallskip

{\bf (iv)} \
$\Omega$ is the singular locus of $F$,
and $w \in \Omega$ \ $\Leftrightarrow$ \ there exist
infinitely many tangent lines to $\Sigma$ passing through $w$.

\end{corollary}

\begin{remark}{}\label{rem:DMS}
Decker, Manolache and Schreyer (cf \cite{DMS} and \ref{prop:E} below) showed
that the partial derivatives of
$f$ define a cubo-cubic Cremona transformation on ${\bf P}(V(14))$,
with base locus $\Omega$ and exceptional
divisor $F$ .
\end{remark}

Furthermore

\begin{proposition}\label{prop:inv}

      $dim \ F = 12$, $deg \ F = 4$, $K_F = {\cal O}_{F}(-10)$.

      $dim \ \Omega = 9$, $deg \ \Omega = 21$.
      $F$ has quadratic singularities along $\Omega\setminus\Sigma$.

      $dim \ \Sigma = 6$, $deg \ \Sigma = 16$,
      $K_{\Sigma} = {\cal O}_{\Sigma}(-4)$.

\end{proposition}

\proof For $F$ it remains to check the statements on singularities.
      Let $f$ be the polynomial defining
      $F$.
      The singularities of $f$ along
      the subscheme defined by the partials are quadratic if and only if
      the subscheme is smooth. But the subscheme defined by the partials
      of $f$ is exactly $\Omega$, which is smooth outside $\Sigma$.
      To compute the invariants of $\Sigma$  we consider the universal exact
sequence of vector bundles on $G=\Gr 36$:
$$0\to U\to V\otimes {\cal O}_{G}\to Q\to 0, $$
\noindent
where $U$ is the universal subbundle.  The $2$-form $\alpha$ restricts
naturally to $U$, i.e. is a section $\alpha_U$ of $(\wedge^2 
U)^{*}\cong\wedge^2U^*$.
The variety $\Sigma$ is therefore nothing but the zero-locus $Z(\alpha_U)$
of this section, so
the class $[\Sigma ]=c_3(\wedge^2U^*)\cap G=(c_1(U^*)c_2(U^*)-c_3(U^*))\cap G$.
   From this description we immediately get that
deg$\Sigma=c_1^6(U^*)\cap\Sigma=16$, and that the canonical divisor
$K_\Sigma= K_G|_\Sigma+c_1(\wedge^2U^*)\cap\Sigma=-4c_1(U^*)\cap \Sigma$.
In particular
$\Sigma$ is a Fano $6$-fold of index $4$.
\qed

\bigskip

\paragraph{}\label{par:de}
We adopt Donagi's notation and
let the {\it pivots} of $w \in {\bf P}^{13} \setminus \Omega$
be the intersection points $\{a,b\}=L_{w}\cap \Sigma$ where $L_{w}$
is the unique secant line through $w$.
Similarly if $w \in F \setminus \Omega$ (= the case when $a = b$),
call $a$ the {\it pivot} of $w$.  When $w\in \Omega$, a {\it pivot} of $w$
is a point $u\in \Sigma$ such that $w$ lies on a tangent line through
$u$.

\paragraph{}\label{par:df}
  From the universal exact sequence it follows that the tangent bundle
$T_{\Sigma}$ is a subbundle of
${\cal H}om(U, U^{*})=U^{*}\otimes U^{*}$.  In fact it is the
subbundle consisting of symmetric
tensors, i.e.
$$T_{\Sigma}={\cal S}ym^2U^{*}.$$
Using the coordinates $[u:X:Y:z]$ above in the point $u = [1,0,0,0]$
on $\Sigma$,
the tangent space $T_u{\Sigma}$ at $u$ to $\Sigma$ is the
span of $\{ [U_{X}]|{\rm rk}X\leq 1\}$, i.e.

$$T_u{\Sigma}={\bf P}^6_u =
<\{ [u:X:0:0] : {}^tX = X \}>$$
The conic in ${\bf P}(U_{0}^\bot)$ defined by the
symmetric matrix $X$ is of course the
conic section $q_{w}(U_{0}^\bot)$ defined in \ref{subsec:rep}.

The following is the $Sp(3)$-analog of
Lemma 3.4 in \cite{D}:

\begin{proposition}\label{prop:tangsp}
      Let $u \in \Sigma$, and let ${\bf P}^6_u$ be, as above,
      the tangent projective space to $\Sigma \subset {\bf P}^{13}$
      at $u$. If $w\in {\bf P}^6_u$, then $q_{w}$ has rank 0, 1, 2 or 3, when
      $w=u$, $w\in \Sigma \setminus u$, $w\in \Omega \setminus\Sigma$ 
or $w\in {\bf
      P}^6_u \setminus\Omega$ respectively.

      \smallskip

      {\bf (i)} \
      $C_{u}:={\Sigma} \cap {\bf P}^6_u$ is a cone
      over the Veronese surface with a vertex $u$;

      \smallskip

      {\bf (ii)} \
      ${\Omega} \cap {\bf P}^6_u$ is a cubic hypersurface.

\end{proposition}
\proof  The tangent cone $C_{u}={\Sigma} \cap {\bf P}^6_u=\{ [U_{X}]|{\rm
rk}X\leq 1\}$. But  ${\rm rk}X, {}^tX = X$ are the equations
of a Veronese surface, hence
$C_{u}={\Sigma} \cap {\bf P}^6_u$ is a cone over a
Veronese surface centered at $u$.  The secant lines to this Veronese
surface fill the
cubic hypersurface in ${\bf P}^6_u$ defined by ${\rm det}X=0$, and
every point on this hypersurface lies on infinitely many secant lines,
so
${\Omega} \cap {\bf P}^6_u$ must coincide with this cubic hypersurface.
\qed

\subsection{Special cycles} \label{cycles}

\paragraph{}\label{para;ea}
We describe some special cycles on $\Sigma$
and on $\check{\Sigma}$.  Again we restrict cycles on $G=\Gr 3V$ to
$\Sigma$.

The restriction of the universal exact sequence on $G$ to $\Sigma$ becomes
$$0\to U\to V\otimes {\cal O}_{\Sigma}\to Q\to 0.$$
The correlation $L_\alpha :V\to V^{*}\quad (v\mapsto\alpha(v,-))$
sets up a natural isomorphism: $Q\cong U^{*}$,
so the universal sequence becomes:
$$0\to U\to V\otimes {\cal O}_{\Sigma}\to U^{*}\to 0.$$
If $\tau_{i}=c_{i}(U^{*})$, then
$$H^{*}(\Sigma)\cong Z[\tau_{1}, \tau_{2},
\tau_{3}]/(\tau_{1}^2-2\tau_{2}, \tau_{2}^2-2\tau_{1}\tau_{3}, \tau_{3}^2 ).$$
Thus the Betti numbers are:
$$1,1,1,2,1,1,1.$$

The classes $\tau_{i}\cap \Sigma$ are represented by cycles that are
restriction of
special Schubert cycles on $G$  to $\Sigma$.
In the notation $\tau_{ijk}=\sigma_{ijk}\cap\Sigma$, we get that
$\tau_{i00}=\tau_{i}\cap \Sigma$.

In projective notation we describe geometrically the cycles  of
Lagrangian planes that
contain a given point, intersect a given line or a given Lagrangian
plane. These belong to the classes $\tau_{300}$, $\tau_{200}$ and
$\tau_{100}$, respectively.
\par
As before, the basic tool will be the correlation
$L_\alpha $.
It is a basic simple fact
that a line in a Grassmannian $\Gr kn$ corresponds to the pencil of $k$-spaces
through a $(k-1)$-space inside a $(k+1)$-space.  Thus a line in $\Sigma$
always correspond to the pencil of $3$-spaces with a common
$2$-subspace in some
$4$-space. From this it follows easily that a plane in $\Sigma$ would
correspond either to the net of $3$-spaces with a common $1$-dimensional
subspace in some $4$-space or to the net of $3$-spaces with a common
$2$-subspace in a $5$-space.
\medskip

\begin{lemma}\label{noplane}

$\Sigma$ contains no planes.

\end{lemma}

\proof Consider the restriction of $\alpha$ to a $4$-space $W$.  The rank of
a skew-symmetric matrix is always even, so the rank of $\alpha_W$ is $0,2$ or
$4$.    Let
$W_\alpha=\cap_{v\in W}{\rm ker}L_\alpha(v).$
    Since $L_\alpha$ is nonsingular, $W_\alpha$ is $2$-dimensional.  Therefore
    dim$W\cap W_\alpha=4-{\rm rank}\alpha_W$ is $0$ or $2$. Let $U$ be a
Lagrangian $3$-subspace of $W$, and let $w\in W \setminus U$. Then
${\rm ker}L_\alpha(v)\cap U= W_\alpha$ so in particular $W_\alpha$ is
contained in $U$.  Thus dim$W\cap W_\alpha=0$ means that there
are no Lagrangian $3$-spaces in $W$, while in case dim$W\cap W_\alpha=2$
there is precisely a pencil.

Similarly, consider the restriction of $\alpha$ to a $5$-space $W$.    The
subspace $W_\alpha=\cap_{v\in W}{\rm ker}L_\alpha(v)$ is $1$-dimensional.
But if it contains a net of Lagrangian $3$-space with a common
$2$-subspace, then
these Lagrangian subspaces fill all of $W$, and hence $W_\alpha$ has dimension
at least $2$, a contradiction.\qed

For $v\in V$ let $V_{v}= {\rm ker}L_\alpha(v)$.

\begin{lemma}\label{points}
For every point $p=<v>\in {\bf P}(V)$
the variety $Q_{p}$ of Lagrangian planes that contain $p$ is a
$3$-dimensional smooth quadric in $\Sigma$.  It is isomorphic to the Lagrangian
Grassmannian $\LG 24$ of Lagrangian subspaces of $V_{v}/<v>$ with
respect to the restriction of $\alpha$ to  $V_{v}/<v>$.
\end{lemma}

\proof  The cycle representing $Q_{p}$ on $\Sigma$ is
$\tau_{3}$ and has degree
$2$.   Any Lagrangian $3$-space that contains $v$ is itself contained
in the $5$-space $V_{v}$. The restriction of $\alpha$ to $V_{v}$ has
kernel $v$, so we may identify $Q_{p}$ with the Lagrangian
Grassmannian with respect to the nondegenerate $2$-form $\alpha_{v}$
induced by $\alpha$ on $V_{v}/<v>$.  This is nothing but a smooth
hyperplane section of a $\Gr 24$, i.e. a smooth quadric $3$-fold. \qed

\medskip

\paragraph{}\label{para;eb} \
The span of the
quadric $Q_{p}$ is a $\Pn {4}\subset {\bf P}(V(14)$ which we denote
by ${\bf P}^4_p$.
Consider the incidence variety

\begin{equation}\label{I_{Q}}
      I_Q=\{([{\bf P}(U)],p)|p\in {\bf P}(U)\}\subset \Sigma\times {\bf P}(V)
\end{equation}

By \ref{points}, this is a quadric bundle over ${\bf
P}(V)$. Now $\Sigma$ spans
${\bf P}(V(14))$, and each ${\bf P}^4_p$ is contained in this span so we may
also consider the incidence

\begin{equation}\label{I_{P}}
I_P=\{(q,p)|q\in {\bf P}^4_p\}\subset{\bf P}(V(14))\times {\bf
P}(V).\label{I_P}
\end{equation}

This is a $\Pn 4$-bundle over ${\bf P}(V)$, which has been studied by
Decker, Manolache and Schreyer.   Its associated rank $5$-bundle is
selfdual, so we describe a construction which is dual to their:
Consider the third exterior power of the Euler sequence on ${\bf P}(V)$
twisted by ${\cal O}_{{\bf P}(V)}(2)$

$$ 0\to \wedge^2 T_{{\bf P}(V)}(-1)
\to\wedge^3 V\otimes {\cal O}_{\Pn 5}(2)
\to \wedge^3 T_{{\bf P}(V)}(-1)\to 0$$
and restrict the contraction $\alpha :\wedge^3 V\to V$ defined by
$$u\wedge v\wedge w\mapsto \alpha (u\wedge v)w+\alpha (v\wedge w)u+\alpha
(w\wedge u)v$$
over ${\bf P}(V)$ to the subbundle
$\wedge^2 T_{{\bf P}(V)}(-1)$.  This restriction, denote it by
$\alpha_U$, is nothing but
the restriction of $\alpha$ over each point $p=<v>$ to $3-$dimensional
subspaces in $V$ which contain $v$.  On the other hand $L_\alpha$ defines a
map
$L_{\alpha_U}: V\otimes
{\cal O}_{{\bf P}(V)}(2)\to {\cal O}_{{\bf P}(V)}(3)$ by $L_{\alpha_U}
(u)(v)=\alpha (u\wedge v)$.
The composition
$$L_{\alpha_U}\circ \alpha_U: \wedge^2 T_{{\bf P}(V)}(-1)\to V\otimes
{\cal O}_{{\bf P}(V)}(2)\to {\cal O}_{{\bf P}(V)}(3)$$ is zero since

$$\alpha (u\wedge v)\alpha (w\wedge u)+\alpha (v\wedge w)\alpha
(u\wedge u)+\alpha
(w\wedge u)\alpha (v\wedge u)=0.$$
The kernel ${\rm ker}L_{\alpha_U}$ is the rank $5$
bundle $\Omega_{{\bf P}(V)} (3)$, so $\alpha_U$ defines a bundle map
$$\alpha_U:\wedge^2T_{{\bf P}(V)}(-1)\to \Omega_{{\bf P}(V)} (3).$$
It is easy to check that this map is surjective as soon as $\alpha $
is nondegenerate. Denote by $E$ the rank $5$ kernel bundle  ker$\alpha_{U}$.
If $U$ is a Lagrangian $3$-space that contains $v$, then $\wedge^3U$
clearly is contained in $E_{p}$ over the point
$p=<v>$.  Thus ${\bf P}(E_{p})$ coincides with the fiber of the
incidence $I_{P}$ over $p$.
\medskip

\begin{proposition}\label{prop:E}(\cite{DMS} Propositions $(1.2)$ and $(1.3)$)
Let $E$ be the rank $5$ kernel bundle of the natural sujective map

$$\alpha_{U}: \wedge^2T_{{\bf P}(V)}(-1)\to \Omega_{{\bf P}(V)}
(3),$$ as above.

Then $E$ has Chern polynomial is $c_t(E)=1+5t+12t^2+16t^3+8t^4$, and
$H^0({\bf P}(V), E)\cong V(14)^{*}$.  Furthermore $E$ is the rank $5$
bundle associated to the $\Pn 4$-bundle $I_{P}$ over ${\bf P}(V)$, and
the projection of $I_{P}$ into the first factor ${\bf P}(V(14))$ is $\Omega$.
\end{proposition}

\medskip

\proof  It remains to compute the Chern polynomial, but this is
straightforward from the construction. The twisted bundle $E(-1)$
coincides with the dual of the bundle
${\cal B}$ defined in
\cite{DMS}.  They show that $B$ is selfdual, so the invariants of the
bundle also follows from their results.\qed

\medskip

\paragraph{}\label{para;ed}
The lines $l\subset \Pn 5$ fall into two cases w.r.t.
$\alpha$. The cycle $\tau_{330}$ representing planes containing a line $l$ is
empty for the general line, while it is a line for each line $l$ in a
Lagrangian plane.  The latter lines are of course precisely the
isotropic lines.

\medskip

\begin{lemma}\label{lines}

If $l$ is a non-isotropic line, then the
variety $\Sigma_{l}$ of Lagrangian
planes that intersect $l$ is a smooth $4$-dimensional variety of degree $8$
which is a quadric bundle inside a rational normal scroll of degree 5, the
Segre embedding of $\Pn 1\times\Pn 4$ in a $\Pn 9$.
If $l$ is an isotropic line the variety $\Sigma_{l}$
    is a singular $4$-dimensional variety of
degree $8$ which is a quadric bundle inside a rational normal scroll
of degree 5
with vertex a line that spans a $\Pn 9$.

\end{lemma}

\proof The variety $\Sigma_{l}$ is represented by the cycle
$\tau_{200}$, the restriction $\sigma_{200}\cap \Sigma$.
Consider the subvariety  $T_{l}\subset\Sigma_{l}$ parameterizing
Lagrangian planes that contain $l$.
The subvariety $T_{l}$ is represented by the class
$\tau_{330}=\tau_{3}^2$ which
is $0$ in $H^{*}(\Sigma)$.  So
there are two cases:

\smallskip

-- the case when $T_{l}$ is positive dimensional (and nonempty), and

\smallskip

-- the case when $T_{l}$ is empty.

\smallskip

In the above coordinates
the first case occurs when
$x_{11}, x_{12}, y_{33}, y_{23}$
vanish.  The equations then reduces to the $2\times
2$-minors of the matrix

$$\pmatrix {u&x_{13}&x_{23}&-x_{22}&y_{13}\cr
           x_{13}&-y_{22}&y_{12}&y_{13}&z\cr}
$$
which define a rational normal $5$-fold scroll with
vertex a line $L(l)$, and

\medskip

\centerline{
$x_{22}x_{33}-x_{23}^2-uy_{11}, \
zx_{33}-y_{11}y_{22}+y_{12}^2, \
x_{13}y_{11}+x_{23}y_{12}+x_{33}y_{13}$.
}

\medskip

\noindent
These quadrics also vanish on $L(l)$.   This is a $4$-dimensional variety of
degree 8 in a $\Pn 9$.   The line $L(l)$ clearly parameterizes
the Lagrangian planes that contain $l$.

The second case occurs when the coordinates
$x_{12}, x_{13}, y_{12}, y_{13}$ vanish.  The equations
then reduces to the $2\times 2$ minors of the matrix

$$\left(\matrix {u&x_{33}&x_{22}&x_{23}&y_{11}\cr
           x_{11}&y_{22}&y_{33}&-y_{23}&z\cr}\right)
$$
which define a smooth rational normal $5$-fold scroll,
the Segre embedding of $\Pn 1\times\Pn 4$, and

$$
x_{22}x_{33}-x_{23}^2-uy_{11},
zx_{11}-y_{22}y_{33}+y_{23}^2,
x_{11}y_{11}-x_{22}y_{22}-x_{23}y_{23}
$$
This is a $4$-dimensional variety of degree $8$ in a $\Pn 9$.
    A bundle over $\Pn 1$ of smooth $3$-dimensional quadrics.
\qed

\medskip

\begin{corollary}\label{Eline}

The restriction $E_l$ of the vector bundle $E$ to a
line $l$ is $E_l=5{\cal O}_l(1)$ when $l$ is non-Lagrangian;
and $E_l=2{\cal O}_l\oplus {\cal O}_l(1)\oplus 2{\cal O}_l(2)$
when $l$ is Lagrangian.

\end{corollary}

\medskip

\proof The restriction is determined by the $2\times 5$ matrices above.\qed

\medskip

We turn to planes.
The isomorphism
$L_{\alpha}:\wedge^3V \cong \wedge^3V^{*} $ induces an isomorphism $\Gr
3V\to\Gr V3$, which composed with the natural isomorphism $\Gr V3\to
\Gr 3V\quad \omega\mapsto {\rm ker}\omega$, defines an involution
$$\iota_{\alpha}:\Gr 3V\to\Gr 3V.$$
The fixed point locus of this involution is of course $\Sigma$, the
set of Lagrangian planes.  On the other hand every involutive pair
$\{P,\iota_{\alpha}(P)\}$ of
planes have a common point
$P\cap\iota_{\alpha}(P)={\rm ker}(\alpha|_{P})$ and is of course contained
in the corresponding correlated hyperplane.

Let $P\subset {\bf P}(V)$ be any plane.  Then the restriction
$\tau_{1}(P)$ of the Schubert cycle $\sigma_{1}(P)$ in $\Gr 3V$ is
clearly a hyperplane
section of $\Sigma$. Thus we have shown
\begin{proposition}\label{plane} If  $P\subset {\bf P}(V)$ is a
non-Lagrangian plane, then
the variety $\Sigma_{P}$ of Lagrangian planes that intersect $P$
coincides with the variety $\Sigma_{\iota_{\alpha}(P)}$ of Lagrangian
planes that intersect $\iota_{\alpha}(P)$ and they define a hyperplane
section of $\Sigma$.
\end{proposition}

\begin{proposition}\label{tangentcone}
For a Lagrangian plane ${\bf P}(U)\subset {\bf P}(V)$ corresponding
to $u\in \Sigma$
the variety $\Sigma_{u}$ of Lagrangian planes that intersect ${\bf P}(U)$
is the hyperplane section defined by the hyperplane
$H_{L_{\alpha}(u)}$,
while the variety of Lagrangian planes that intersect ${\bf P}(U)$ in a
line is a cone $C_{u}$ over a Veronese surface of degree $4$ with
vertex at $u$ in a
$\Pn {6}=\Pn {6}_{u}$, the projective tangent space to $\Sigma$ at $u$.
\end{proposition}

\proof  $\Sigma_{u}$ is a hyperplane section with ${\bf P}(U)$
invariant under the involution $\iota_{\alpha}$.   The lines in
$\Sigma$ through $u$ represent precisely the pencils of Lagrangian
planes that contain ${\bf P}(U)$.  On the other hand these lines
generate precisely the tangent cone $C_{u}$ to $\Sigma$ at $u$, so the lemma
follows from \ref{prop:tangsp}.  \qed

   Recall the decomposition $V=U_{0}\oplus U_{1}$ from
\ref{subsec:rep}.  The restriction of the global sections
$V(14)^{*}$
of $E$ to the Lagrangian plane $P={\bf P}(U_{0})$ decomposes as the
restriction of the decomposition

$$ \wedge^3{V^{*}}=
\wedge^3U_{1}^{\bot}\oplus(\wedge^2U_{1}^{\bot}\otimes U_{0}^{\bot})\oplus
(U_{1}^{\bot}\otimes\wedge^2U_{0}^{\bot})\oplus \wedge^3U_{0}^{\bot}$$
to $V(14)^{*}$.  By the natural map $q(U_{0})$ this decomposition
becomes
$$H^0(P,E|_P)\cong \wedge^3U_{1}^{\bot}\oplus U'\oplus
{\cal S}ym^2U_{0}^{*},$$
where the first summand consists of the constant forms, the second is the
restriction of the summand $(\wedge^2U_{1}^{\bot}\otimes
U_{0}^{\bot})$, while the last summand is the quadratic forms on
$P$.  The restriction of the vector bundle $E$ to $P$ therefore
decomposes as a sum of two line bundles and a rank $3$ bundle.
Thus

\begin{proposition}\label{Eplane}

For a Lagrangian plane $P= {\bf P}(U)\subset {\bf P}(V)$
the restriction of the vector bundle $E$ to $P$ is
\ $E|_P={\cal O}_P\oplus {\cal O}_P(2)\oplus E_P$ \
where $E_P$ is a rank 3 vector bundle with Chern polynomial
$c_t(E_P)=1+3t+6t^2$.
\end{proposition}
\proof  The Chern polynomial follows by a direct calculation.
\qed

\subsection{Singular hyperplane sections}\label{hyperplane sections}

\paragraph{}\label{para;fa}

In this section we describe the hyperplane sections of $\Sigma$
corresponding to the different $Sp(3,{\bf C})$-orbits in ${\bf
P}(V(14)^{*})$.
In particular we describe their  singular locus.  Four maps are
crucial.
First there is the basic correlation
$$L:  \Sigma\to \check{\Sigma},\quad {\bf P}(U)\mapsto {\bf
P}(L_{\alpha}(U)).$$
The second map is the involution
$$\iota_{\alpha}:\Gr 3V\to \Gr 3V.$$

The third map, which we call the vertex map, is
$$v:\Omega \setminus{\Sigma}\to {\bf P}(V)\quad q\mapsto
\pi_{2}\pi_{1}^{-1}(q),$$
where $\pi_{1}$ and $\pi_{2}$ are the projections from the
incidence $I_{P}$.  The point $v(w)\in {\bf P}(V)$, is called the
vertex of $w\in \Omega$.
The third map is the pivot map:
$$piv:F \setminus{\Omega}\to \Sigma\quad p\mapsto
{\rm the\; (unique)\; pivot\; of\; }p.$$
The corresponding maps on the dual space
are marked with a $\ast$.
Notice that we have the following relations:
$$L^{-1}\circ piv^{*}=piv\circ L^{-1}:\check{F} \setminus
\check{\Omega}\to \Sigma$$
and
$$L^{-1}\circ v^{*}=v\circ
L^{-1}:\check{\Omega} \setminus
\check{\Sigma}\to {\bf P}(V).$$
The first of these will be denoted by $u$:
$$u=L^{-1}\circ piv^{*}:\check{F} \setminus
\check{\Omega}\to \Sigma$$
while, by abuse of notation, the second one will be denoted by
$v$:
$$v=L^{-1}\circ v^{*}: \check{\Omega} \setminus
\check{\Sigma}\to {\bf P}(V).$$
For the point $\omega \in {\bf P}(V(14)^{*}$, denote by
${\bf P}^{12}_\omega \subset {\bf P}(V(14))$ the hyperplane
defined by $\omega$, and denote by $H_\omega = {\bf P}^{12}_\omega \cap \Sigma$
the corresponding hyperplane section of $\Sigma$.

\begin{lemma}\label{dualsing}
The hyperplane section
$H_{\omega}={\bf P}^{12}_{\omega}\cap \Sigma$
is singular precisely in the points
$u\in \Sigma$ such that the tangent space $\Pn 6_{L(u)}$ to
$\check{\Sigma}$
contains $\omega$.
\end{lemma}

\proof The hyperplane section $H_{\omega}$ is singular at $u$ if and only if
it contains the tangent space $\Pn 6_{u}$, which is equivalent by
$L_{\alpha}$ to $\omega\in \Pn 6_{L(u)}$.\qed

\begin{proposition}\label{singh}
{\bf (i)} \
When $\omega\in {\bf
P}(V(14)^{*})  \setminus \check{F}$,
then $H_{\omega}={\bf P}^{12}_{\omega}\cap \Sigma$ is  smooth.

{\bf (ii)} \
When $\omega\in \check{F} \setminus
\check{\Omega}$, then $H_{\omega}$ is  singular precisely in the point
$u(\omega)$.

{\bf (iii)} \
When
$\omega\in \check{\Omega} \setminus \check{\Sigma}$,
then $H_{\omega}=\Sigma_{P_{1}}=\Sigma_{P_{2}}$ for an involutive pair
of planes  $P_1(\omega)$ and
$P_2(\omega)=\iota_{\alpha}(P_{1}(\omega))$. Furthermore $H_{\omega}$
is singular along a smooth
quadric surface $Q_{\omega}$ in $P_{v(\omega)}$ which parameterizes
the set of  Lagrangian planes passing through $v(\omega)$ and
intersecting $P_1$  and $P_2$ in a line.

{\bf (iv)} \
When $\omega\in \check{\Sigma}$, then $H_{\omega}$ is
singular along $C_{u}$ the cone over a Veronese surface with vertex at
$u=L^{-1}(\omega)$.
\end{proposition}

\proof We check one point $\omega $ in each orbit and use
\ref{dualsing} to find the singular points of $H_{\omega}$ as the
points $u\in \Sigma$ such that the tangent spaces $\Pn 6_{L(u)}$
contain $\omega$.
   Clearly $H_{\omega}$ is smooth
    \ {\it iff} \ $w\in \Pnd {13} \setminus \check{F}$,
    while $H_\omega$ is singular precisely at
    $u(\omega)$.  On the other hand \ref{prop:tangsp} implies that
$H_{\omega}$ is singular along the cone $C_{u}$ with vertex at
$u=L^{-1}(\omega)$ when $\omega\in \check{\Sigma}$.

It remains to check a point $\omega\in\check{\Omega} \setminus \check{\Sigma}$,
and we  start with a point $p$ such that $p=v(\omega)$.
In the usual coordinates we consider the point $p=<e_1>\in {\bf
P}(V)$.  Any Lagrangian plane
through $p$ lies in the hyperplane $x_{4}=0$.  We construct the
singular locus of a hyperplane section $H_{\omega}$ with $\omega\in
{\bf P}^4_{L_{\alpha}(p)}\setminus Q_{L_{\alpha}(p)}$.   Consider the 
two planes

$$P_1 = {\bf P}(<e_1,e_2,e_5>)\; {\rm and}\; P_2 =
\iota_{\alpha}(P_{1}) = {\bf
P}(<e_1,e_3,e_6>)$$

\noindent
passing through $p$ and spanning this hyperplane.
For any $(t) = (t_2:t_5)$ and for any $(s) = (s_3:s_6)$
the plane
${\bf P}^2(s;t) ={\bf P}( <e_1,t_2e_2+t_5e_5,s_3e_3+s_6e_6>)$
is Lagrangian and intersect both $P_1$ and $P_2$ in a line.
So the corresponding point
$$u(s;t) =< e_1 \wedge (t_2e_2+t_5e_5) \wedge (s_3e_3+s_6e_6) >\in
\Sigma.$$

Consider $P_{u(s;t)}$, the
projective tangent space to $\Sigma$ at $u(s;t)$,
$(s;t) \in {\bf P}^1 \times {\bf P}^1$.
For a 3-vector $u = x \wedge y \wedge z \in \Sigma$
the projective tangent space
$P_u$  is spanned by

$$\{w\in\wedge^2 U\wedge V|\alpha (w)=0\}$$
where $U=<x,y,z>$.
Applied to the tangent spaces
${\bf P}^6(s;t)$, we find that the point $w=e_{125} + e_{163}$ is
a common point.  For this point only the coordinate $x_{23}$
is nonzero, so by the equations, it does not lie on $\Sigma$.
On the other hand all the $u(s;t)$ are ``pivots''
of the point $w$, i.e. for any $u(s;t)$ the line
$<u(s;t),w)>$ is tangent to $\Sigma$.
Thus the point $w$ will belong to $\Omega \setminus \Sigma$. In fact $w\in
{\bf P}^4_{p}$ and every pivot of $w$ is a Lagrangian plane through $p$.
Since $Q_{p}=\Sigma\cap {\bf P}^4_{p}$  is a quadric 3-fold, the set of pivots
of $w$ is the intersection of $Q_{p}$ and the polar of $w$ w.r.t.
$Q_{p}$, i.e. a quadric surface.
Therefore $S_w = \{ u(s;t): (s;t) \in {\bf P}^1 \times {\bf P}^1 \}$
is precisely the set of pivots of $w$.

Let $\omega=L(w)=x_{452}+x_{436}\in \check{\Omega} \setminus \check{\Sigma}$,
then by \ref{dualsing} the hyperplane section $H_{\omega}$ is
singular precisely along
$S_{L^{-1}(\omega)}=S_{w}$. Finally, any Lagrangian subspace
$U'\subset V$ that
intersects $<e_{1},e_{2},e_{5}>$, intersects also $<e_{1},e_{3},e_{6}>$,
so it is of the form
$$<ae_{1}+be_{2}+ce_{5},a'e_{1}+b'e_{3}+c'e_{6}, v>$$
so $\omega(\wedge^3U')=0$, i.e. the corresponding point $u'\in
\Sigma_{P_{1}}$.
\qed
We give a more precise description of the type of singularities of
hypersections.

\begin{proposition}\label{singha}

{\bf (i)} \
If $\omega\in \check{F} \setminus \check{\Omega}$
then the hyperplane section $H_\omega \subset \Sigma$ has a quadratic
singularity at the pivot $u = u(\omega) = L(piv^{*}(\omega)) \in \Sigma$.

{\bf (ii)} \
If $\omega \in \check{\Omega} \setminus \check{\Sigma}$
then the hyperplanes section $H_{\omega} \subset \Sigma$
has a quadratic singularity
{\it along} the quadric surface $Q_{\omega} \subset \Sigma$.

{\bf (iii)} \
Let $\omega \in \check{\Sigma}$, and let $w = L^{-1}(\omega) \in \Sigma$.
Then the projective tangent cone
is a cone $C_{w}$ over a Veronese surface with vertex $w$.
The hyperplane section $H_{\omega} \subset \Sigma$ has
a double singularity of rank $3$ along the punctured cone
$C_w \setminus \{ w \}$.
\end{proposition}

\medskip

\proof  We choose one representative for each orbit. Let $u = [1,0,0,0]$.
By \ref{subsec:rep}, symmetric $3\times
3$-matrices $X$ and $Y$ such that
$[1,X,Y,z] \in {\bf C}^{13} \subset {\bf P}(V(14)) = <\Sigma>$
form local coordinates at the point $u = [1,0,0,0]$
in which a neighborhood ${\Sigma}^o \subset \Sigma$
has a local parameterization
$[1,X,{\wedge}^2 \ X, det(X)]$.
The tangent space $T{\Sigma}|_u$ is parameterized by the
linear matrix
$X = (x_{ij})$,
and the projection $[1,X,Y,z] \rightarrow X$
sends the neighborhood ${\Sigma}^o$
to ${\bf C}^6(x_{ij}) \cong T{\Sigma}|_u$.

{\bf (i)} \
We choose $\omega\in T{\Sigma}|_u^{\bot}$ such that
$H_{\omega}=\{y_{11}+y_{22}+y_{33}=0\}$.
Clearly this hyperplane contains the projective tangent space
${\bf P}^6_u = {\bf P}(T{\Sigma}|_u)$ to $\Sigma$ at $u$,
in particular $u \in H_\omega$.
In matrix coordinates $\omega$ is represented by the matrix
$$
Y =
\left(
\begin{array}{ccc}
1 & 0 & 0 \\
0 & 1 & 0 \\
0 & 0 & 1
\end{array}
\right).
$$
Since $rank(Y) = 3$,
form $\omega \in \check{F}\setminus\check{\Omega}$.
Since on ${\Sigma}^o$ the parameters
$y_{ij}$ are  in fact the minors of the matrix $X$,
the isomorphic projection of the neighborhood
$H_{\omega}^o = {\Sigma}^o \cap {\bf P}^{12}_\omega\subset H_\omega$
of $u$ on $T{\Sigma}|_u = {\bf C}^6(x_{ij})$
is a hypersurface with a local equation
$y_{11}+ y_{22}+ y_{33} =
Q(X)
=
x_{22}x_{33}-x_{23}^2
+
x_{11}x_{33}-x_{13}^2
+
x_{11}x_{22}-x_{12}^2
= 0$.
Therefore $mult_u \ H_{\omega} = deg \ Q(X) = 2$.
Moreover $u$ is a quadratic singularity on
$H_\omega$ since the quadric $\{Q(X) = 0\}\subset {\bf P}(T{\Sigma}|_u)$
has maximal rank $= 6={\rm dim}H_{\omega}+1$.

\medskip

{\bf (ii)} \
Next, we choose $\omega\in\check{\Omega} \setminus \check{\Sigma}$ such that
$H_{\omega} = (y_{11}+y_{22} = 0)$.
Thus $u\in Q_{\omega}$, the singular locus of $H_{\omega}$.
Now, the equation of the projectivized tangent cone
to $H_{\omega}$ at $u$ becomes:
$y_{11}+ y_{22} =
Q(X)
=
x_{22}x_{33}-x_{23}^2
+
x_{11}x_{33}-x_{13}^2
= 0$.
Therefore $mult_u \ H_{\omega} = deg \ Q(X) = 2$,
but in this case $rank(Q) = 4 < 6$.  On the other hand $H_{\omega}$
is singular along the
surface $Q_{\omega}$, so the singularity has
maximal rank, in particular it is quadratic along
$Q_{\omega}$.

\medskip

{\bf (iii)}  \
Let $\omega \in \check{\Sigma}$ such that $H_{\omega} = (y_{11} =
0)$, and let $w = L^{-1}(\omega)$.
The tangent cone at $w$ of $H_{\omega}$ is a cone
$C_w$ over a Veronese surface
with vertex $w$.  The point $u = [1,0,0,0]$ lies on this cone and is
different from $w$.
The local equation of
$H_{\omega}$ at $u$ is:
$y_{11} =
Q(X) = x_{22}x_{33}-x_{23}^2
= 0$.
Since $rank(Q) = 3$ the hyperplane section
$H_{\omega}$ has a quadratic singularity
along the $3$-dimensional punctured cone
$C^*_w = C_w \setminus \{ w \}$.
\qed

\medskip

\paragraph{}\label{para;fb}
Now, we turn to the involutive pair of planes
$P_1(\omega),P_2(\omega)$ that appear in \ref{singh}.
It follows from \ref{singh} (iii)  that the union of the planes
$P_1(\omega)\cup P_2(\omega)$ is
the set of points $q$ in ${\bf P}(V)$, such that
${\bf P}^4_q\subset {\bf P}^{12}_\omega$
(or, equivalently, $Q_q \subset H_\omega$).
This fact fits in the description of the incidence

\begin{equation}\label{J}
J=\{(p,w)\in {\bf P}(V)\times {\bf P}(V(14)^{*})|{\bf P}^4_p\subset
{\bf P}^{12}_\omega\}
\subset {\bf P}(V)\times {\bf P}(V(14)^{*}).
\end{equation}

Let ${\pi}:J \rightarrow \Pn 5={\bf P}(V)$
and ${\psi}:J \rightarrow \Pnd {13}={\bf P}(V(14)^{*})$ be the two
projections of $J$.

\medskip

\begin{proposition}\label{incidence}
The image of the second
projection is precisely the quartic
$\check{F} \subset \check{\bf P}^{13}$.
The fiber of the second projection over a point $\omega\in
\check{F}\setminus\check\Omega$
is the smooth conic
section $q_\omega(U)\subset{\bf P}(U)$, where ${\bf P}(U)$ is the plane
of the pivot $u(\omega)$.
The fiber of $\psi$ over a point $\omega$ on
$\check{\Omega}\setminus\check{\Sigma}$ is the union of the two
planes $P_1(\omega)\cup
P_2(\omega)$.  The fiber of $\psi$ over a point $\omega\in
\check{\Sigma}$ is the Lagrangian plane ${\bf P}(U)$ of
$u=L^{-1}(\omega)$.
\end{proposition}

\medskip

\proof The last statement follows from
\ref{tangentcone}, while the case
$\omega\in\check{\Omega}\setminus\check{\Sigma}$
follows as explained above from \ref{singh} (iii).
The remainder of the proposition we reformulate:

\bigskip

\begin{proposition}\label{inca}
Let $\omega \in {\bf P}(V(14)^{*})$, and let $p \in {\bf P}(V)$.
Then ${\bf P}^4_p\subset{\bf P}^{12}_\omega$
$\Leftrightarrow$ $\omega\in T_{[U]}^{\bot}$ for some Lagrangian subspace
$U$ with $p\in  {\bf P}(U)$ and $q_{\omega}(U)(p)=0$.
\end{proposition}
\proof By abuse of notation we do not distinguish between $\omega\in {\bf
P}(V(14)^{*})$ and any nonzero vector in $V(14)^{*}$ representing it.
Thus we
consider $\omega$ as a section of the vector
bundle $E$, and analyse its restriction to Lagrangian planes.

Let $P={\bf P}(U)\subset{\bf P}(V)$ be a Lagrangian plane and $u\in \Sigma$ the
corresponding point.
According to \ref{Eplane} the restriction of the vector bundle $E$ to
$P$ decomposes into three direct summands:
$$E|_P={\cal O}_P\oplus {\cal O}_P(2)\oplus E_P$$
where $E_P$ is a rank 3 vector bundle with Chern polynomial
$c_t(E_P)=1+3t+6t^2$.  Therefore the restriction $\omega_{P}$ of
$\omega$ to $P$ decomposes into
$\omega_{P}=a\oplus b\oplus c$
where $a$ is a constant, while $b$ is a quadratic form and $c$ is a
section of the rank $3$ bundle $E_{P}$.

\begin{lemma}\label{rplane}
\smallskip
{\bf (i)} \
If $\omega\in \check{\Sigma}$, such that the dual pivot $u(\omega)=u$, then
$a(\omega)=b(\omega)=c(\omega)=0$.

{\bf (ii)} \
$a(\omega)=c(\omega)=0$ if and only if $\omega\in \check{F}$ and the
dual pivot $u(\omega)=u$.

{\bf (iii)} \
$u \in {\bf P}^{12}_\omega$ if and only if  $a=0$.

\end{lemma}

\proof  Let $u\in \Sigma$ with corresponding Lagrangian plane $P={\bf
P}(U)$.  The hyperplane ${\bf P}^{12}_\omega$ contains $\Sigma_u$ (cf.
\ref{tangentcone}) if and only if $a=b=c=0$,
so (i) follows.
The hyperplane ${\bf P}^{12}_\omega$ contains the tangent cone
$C_{u}$ at $u$ (the cone over a
Veronese surface) if and only if the restriction $\omega_{P}\in {\cal
S}ym^2(U)^{*}$, i.e.
$a(\omega_{P})=c(\omega_{P})=0$, so (ii) follows.  Finally
$a=0$ if and only if the ${\bf P}^{12}_\omega$
passes through the vertex $u$ of $\Sigma_u$.\qed

If  $\omega\in \check{F}$ and the dual pivot $u(\omega)=u$, then the
quadratic form $b(\omega_{P})$ is nothing but the quadratic form
$q_{\omega}(U)$.   Let $p\in P={\bf P}(U)$.  Then
$\omega(p)=0$ if and only if $\omega\in \check{F}$ with dual pivot
$u(\omega)=u$ and $q_{\omega}(U)(p)=0$ by \ref{rplane} (ii).  On the
other hand $\omega(p)=0$ if and only if
   ${\bf P}^4_p\subset{\bf
P}^{12}_\omega$ so the proposition follows.\qed

\subsection{Linear sections, $Sp(3)$-dual sections and vertex
varieties}\label{linear sections}

In this section we explore the incidence correspondences $I_P$, $I_Q$ and
$J$  of the previous two sections to study the relations between
linear sections of $\Sigma$, the
dual linear sections of $\check{F}$ and certain vertex varieties.

\paragraph{}\label{para:ga}

For $2\leq k\leq 5$, let
${\bf P}^{13-k}
\subset {\bf P}(V(14)$ be a general linear subspace of codimension
$k$,  and let
${\Pi}^{k-1} = ({\bf P}^{13-k})^{\perp}\subset {\bf P}(V(14)^{*}) $
be the orthogonal subspace of hyperplanes
which pass through ${\bf P}^{13-k}$.
Let $$X^{6-k} = \Sigma \cap {\bf P}^{13-k},\quad
\check{F}(X^{6-k})={\Pi}^{k-1}\cap \check{F}, $$
and let $${\check{\Omega}}(X^{6-k}) = \check{F}(X^{6-k}) \cap \check{\Omega}
= {\Pi}^{k-1} \cap \check{\Omega}.$$ Thus the superscript always
indicates the dimension of the linear section $X$.
We call $\check{F}(X^{6-k})$ the {\it $Sp(3)$-dual} section to $X^{6-k}$.

We restrict our attention to general linear subspaces, more precisely
we will assume that $X^{6-k}$ is smooth.
For singular sections $X^{6-k}$, there are obviously similar results.

\medskip

\begin{lemma}\label{singP}
Let $P\subset {\bf P}(V(14))$ be a linear subspace and let
$\omega \in P^{\perp}\cap\check{F}$
such that a dual pivot $u\in \Sigma$ of $\omega$
lies in $P$.
Then $u \in Sing \ P\cap\Sigma$.

\end{lemma}

\proof Let $H_\omega = \Sigma \cap {\bf P}^{12}_\omega$ be the hyperplane
section of $\Sigma$ defined by $\omega$.
By \ref{dualsing}, $H_\omega$ is singular at the pivot $u$.
Since $\omega \in P^{\perp}$,
the variety $P\cap \Sigma \subset H_\omega$ is a
complete intersection of $H_\omega$ and hyperplanes
which pass through the singular point $u$.
Therefore $P\cap \Sigma$ is singular at $u$.
\qed

\medskip

It follows from \ref{singP} that if
$X^{6-k} = \Sigma \cap {\bf P}^{13-k}$
is smooth then $u(\omega) \not\in {\bf P}^{13-k}$
for any $\omega \in \check{F}(X^{6-k})\setminus\check{\Omega}(X^{6-k})$.
Combined with \ref{prop:inv} this yields:
\medskip

\begin{proposition}\label{subcanonical}

Let $X^{6-k}$ be a smooth $(6-k)$-dimensional linear section of $\Sigma$.

If $2 \le k \le 4$, then $\check{\Omega}(X^{6-k}) = \emptyset$ and
$\check{F}(X^{6-k}) \subset {\Pi}^{k-1}$ is a smooth quartic $(k-2)$-fold
(e.g. when $k = 2$ then  $\check{F}(X^{6-k})$
is a set of $4$ points each with multiplicity $1$);
and if $k = 5$ then
$Sing \ \check{F}(X^{1}) = {\check{\Omega}(X^{1})} = \{
{\omega}_{1},...,{\omega}_{21} \}$
is a set of $21$ ordinary double points (nodes)
of the quartic threefold $\check{F}(X^{1})$.

The linear section $X^{6-k}$ is subcanonical, more precisely
   $-K_X = (-4+k)H$
where $H$ is the class of the hyperplane section.
When
$k\leq 4$ and $X^{6-k}$ is general, then
${\bf Pic}(X^{6-k}) = {\bf Z}[H]$
\end{proposition}

\proof Only the last statement remains to be shown, but this follows
from \cite{Mo}  (cf.
also \cite{Mu1}).
\qed

\paragraph{Two quadric bundle fibrations}\label{gb}

Let $1\le k\le 6$, let $P=\Pn {13-k}$ be a general linear space in
${\bf P}(V(14)$
as above and consider the incidence variety

$$I_{P}(k) =
\{(x,y)\in X^{6-k}\times{\bf P}(V)|{\rm dim}({\bf P}^4_y\cap\Pn {13-k})\geq
k-1, x\in
Q_y\}\subset I_P$$

We denote the image of the second projection in ${\bf P}(V)$ by $Y(k)$ and
call it the vertex variety of $X^{6-k}$.
The linear space $P=\Pn {13-k}$ is defined by $k$ linear forms.  These linear
forms may be pulled back to global sections of the vector bundle $E$
(cf. \ref{prop:E}), therefore they define a map
$\Phi_{P}: k{\cal O}_{{\bf P}(V)}\to E$  and  $Y(k)$ is just the
$(k-1)$th degeneracy
locus
$$Y(k)= \{y|{\rm rank}\Phi_{P}(y)\leq
k-1\}\subset {\bf P}(V).$$


Thus $Y(k)$ has dimension $k-1$ and degree given by Porteous' formula 
\cite{Fu}:
$${\rm deg}Y(k)=c_{6-k}(E),\qquad k=1,\ldots,5$$

The fibers of the map
$$I_{P}(k) \to Y(k)$$ are all quadrics.  In fact

\begin{proposition}\label{prop:quadric bundle} When $2\leq k\leq 4$
and $\Pn {13-k}$ is a general linear subspace in ${\bf P}(V(14))$, then the map
$I_P(k) \to Y(k)$ is a quadric bundle map, i.e. the fibers are all
quadrics of dimension $k-2$.\end{proposition}
\proof When the $k$ sections of $E$ are general, ${\rm rank}\Phi_{P}(y)\geq
k-1\}$ for every $y$.  Therefore ${\rm dim}({\bf P}^4_y\cap\Pn {13-k})=
k-1$ for every $y\in Y(k)$ and the proposition follows.\qed

The fibers of the other projection
   $$I_P(k)\to X^{6-k}$$
   are the intersections of a Lagrangian
plane of $X^{6-k}$ with the subvariety $Y(k)$.  This intersection is the
rank $k-1$ degeneracy locus of the vector
bundle map $\Phi$ restricted to these Lagrangian planes.

\paragraph{}\label{gf}

Consider now $k$ sections of $E$, and let $\Pn {k-1}$ be the linear
subspace of $\Pnd {13}$ spanned by these sections.  Let $M(k)$ denote the
(locally) $k\times 5$ matrix with rows defined by the
sections as above.

\medskip

\begin{lemma}\label{lagplanes}

Let $u\in (\Pn {k-1})^{\bot}\cap \Sigma$ and let $P$ be the
corresponding Lagrangian plane in ${\bf P}(V)$, then
$Y(k)$ contains $P$ if $k=5$, it intersects  $P$ in
a quintic curve if $k=4$ and in a finite scheme of length $12$ if
$k=3$.\medskip
Let  $\omega\in \check{F}\cap \Pn {k-1}$, and let $P_{u}$ be the
Lagrangian plane of $u=u(\omega)$ the
dual pivot. If $u\notin (\Pn {k-1})^{\bot}\cap \Sigma$, then $P_{u}\cap
Y(k)=\{q_{\omega}=0\}$, the conic section
defined by $\omega$.
\end{lemma}

\proof The restriction to $P$ of
$M(k)=(a(\omega);b(\omega);c(\omega))$ has,  by
\ref{rplane} (iii), $a(\omega)=0$.  Therefore the rank of the matrix is
at most $(k-1)$ along a subvariety of class $c_{5-k}(E)\cap P$.
If $\omega\in \check{F}\cap \Pn {k-1}$, then the rank of the matrix is
at most $(k-1)$ along the conic section $q_{\omega}$. \qed

\begin{proposition}\label{fiber2} The fibers of the projection
   $$I_P(k)\to X^{6-k}$$
   are Lagrangian planes if $k=5$, plane quintic curves if $k=4$, and
   finite planar subschemes of length $12$ if $k=3$.
\end{proposition}

\paragraph{}\label{gg}

We turn to the second incidence variety $J$ of (\ref{J}).  It involves the
$Sp(3)-$dual  $\check{F}(X^{6-k})=\check{F}\cap
(P^{13-k})^{\bot}$.  Consider first

$$IJ_{P}(k) =
\{(x,y,\omega)\in X^{6-k}\times\Pn 5\times \check{F}(X^{6-k})|{\rm
dim}(\Pn 4_y\cap\Pn {13-k})\geq
k-1, x\in Q_y, {\bf P}^4_y\subset {\bf P}^{12}_\omega\}$$
and its projection onto the second and third factor
$$J_{P}(k)=\{(y,\omega)\in \Pn 5\times \check{F}(X^{6-k})|{\rm dim}(\Pn
4_y\cap\Pn {13-k})\geq
k-1,  {\bf P}^4_y\subset {\bf P}^{12}_\omega\}.$$

\medskip

\begin{lemma}\label{isoproj}

When $k\leq 4$ and $({\bf P}^{13-k})^{\bot}$ does not
intersect $\check{\Omega}$, then the projections
$IJ_{P}(k)\to I_{P}(k)$ and
$J_{P}(k)\to Y(k)$ are isomorphisms.
\end{lemma}

\proof When $({\bf P}^{13-k})^{\bot}$ does not intersect $\check{\Omega}$, then
for every $y\in Y(k)$ there is a unique $\omega\in
\check{F}(X^{6-k})$ such that
${\bf P}^4_{y}\subset {\bf P}^{12}_{\omega}$.  \qed

\begin{proposition}\label{conic bundle}
When $\check{F}(X^{6-k})\subset \check{F} \setminus
\check{\Omega}$, then $Y(k)$ has a conic
bundle structure over $\check{F}(X^{6-k})$, where every conic is smooth.
\end{proposition}

\proof  Since $Y(k)\cong J_{P}(k)$ by \ref{isoproj}, it suffices to
consider the map $J_{P}(k)\to \check{F}(X^{6-k})$.  Let $\omega \in
\check{F}(X^{6-k})$.
Since $\omega\in \check{F} \setminus\check{\Omega}$, there is a double pivot
$u=u(\omega)\in
\Sigma$.
According to \ref{incidence} (or \ref{lagplanes}) the Lagrangian plane $P(U)$
corresponding to $u$ intersects
$Y(k)$ in a smooth conic section $\{q_\omega(U)=0\}$.\qed

The invariants of the vertex variety $Y(k)$ is easily computed from
the degree and
the conic bundle structure.
The only involved case is $k=4$, i.e. when $Y(4)$ is a threefold
conic bundle over the $K3$-surface $\check{F}(X^{2})$.  These
threefolds form one of only two families of smooth threefold conic
bundles in $\Pn
5$ and is extensively studied in \cite{BOSS}.

\begin{corollary}\label{s=9} (\cite{BOSS})
For a general linear surface section $X^2\subset \Sigma$
the vertex variety $Y(4)$ is a threefold smooth conic bundle over a
quartic surface.
It has degree 12 and
intersection
numbers:
     $$H^{3-i}\cdot (H+K)^i=c_{2+i}(I_{Y^s}(5h)).$$
\end{corollary}

\medskip

\begin{corollary}\label{s=10}
For a general $3$-fold linear section $X^3\subset\Sigma$
the vertex variety $Y(3)$ is a minimal smooth conic bundle of
degree 16 over a
quartic plane curve.
\end{corollary}

\begin{corollary}\label{s=11}

For a general $4$-fold linear section $X^4\subset \Sigma$
the vertex variety $Y(2)$ is 4 disjoint conic sections.
\end{corollary}

\paragraph{}\label{para:gh}
We  summarize these results in the following table:

\medskip
\begin{center}
\begin{tabular}{|c|c|c|c|c|c|} \hline
$k $ & $X^{6-k}$&$Y(k)$&$\check{F}(X^{6-k})$\\ \hline\hline
$1$&Fano 5-fold of index 3& $\emptyset$&$\emptyset $ \\\hline
$2$&Fano 4-fold of index 2&4 disjoint conics& 4 points \\\hline
$3$&Fano 3-fold of index 1&surface of degree 16&plane quartic curve \\\hline
$4$&K3-surface&3-fold of degree 12&quartic surface  \\\hline
$5$&Canonical curve of genus 9&quintic hypersurface&21-nodal quartic
3-fold  \\\hline

\end{tabular}
\end{center}

\bigskip

\subsection{The Fano 3-fold linear sections}\label{fano threefold}

\bigskip

The general $3$-dimensional linear section $X=X^3$ of $\Sigma$ is a
prime Fano 3-fold of genus $9$
(see \cite{M}, \cite{M2}). The vertex variety $Y=Y(3)$ of $X$ is a
conic bundle over the $Sp(3)$-dual plane quartic curve $\check{F}(X)$
by \ref{conic bundle}.
Below we give several relations between
the Fano 3-fold $X$ and its vertex surface $Y$.  For a different 
approach to these relations see \cite{I}.


\begin{proposition}\label{lines10}
      For the general smooth $3$-fold linear section $X$ of 
      $\Sigma$ there is a one to
    one correspondence between the set of lines in $X$, the $5$-secant
    lines to its vertex variety $Y$ and sections $C$ with minimal self
intersection $C^2=3$ of $Y$ as a
    conic bundle.\end{proposition}

\proof A line $L$ in $\Sigma$  parameterizes the pencil of Lagrangian planes
through an isotropic line $l$, in the notation of \ref{lines}
$L=L(l)$. The union of these Lagrangian planes is
a $3$-space $\Pn 3_{l}$.   We will show that $l$ is a $5$-secant line
to $Y$ while  $\Pn 3_{l}\cap Y$ contains a section $C_{L}$ of
the conic bundle $Y\to \check{F}(X)$.  First we consider the
restriction of the vector
bundle map $\Phi$ to $l$.
It is defined by a matrix with rows
$$M_l=\pmatrix {a_{1}&a_{2}&b_1&c_{1}&c_{2}}.$$
Then
$$\Phi_l:{\cal
O}_l\to 2{\cal O}_l\oplus {\cal O}_l(1)\oplus 2{\cal O}_l(2).$$
Thus the $a_i$ are constants, while $b_i$ and $c_{i}$ are linear and
quadratic respectively.
A point  $\omega\in\Pnd {13}$ define a section of $E$, i.e.
a matrix $M_l(\omega)$.  In the notation of \ref{lines} we have

\begin{lemma}\label{rline}
If $l$ is an isotropic line,
   then $\Sigma_{l}$ is contained in ${\bf P}^{12}_{\omega}$
if and only if $M_l(\omega)=0$.  The line $L(l)$ parameterizing the
pencil of Lagrangian
planes through $l$ is contained in the hyperplane ${\bf P}^{12}_\omega$
defined by $\omega$,
if and only if $a(\omega)=0$.
\end{lemma}

\proof The first statement is obvious, while $L(l)$ is the vertex of
$\Sigma_{l}$ (cf \ref{lines}), so clearly, ${\bf P}^{12}_\omega$
contains the vertex $L(l)$ of
$\Sigma_l$ precisely when $a=0$.
\qed

When $l$ is a non-isotropic line, then the
row matrix

$$M_l=\pmatrix {a_{1}&a_{2}&a_{3}&a_{4}&a_{5}}$$
has only linear entries.

\begin{lemma}\label{rnline}
If $l$ is a non-isotropic line,
then $\Sigma_{l}$ is contained in ${\bf P}^{12}_{\omega}$
if and only if $M_l(\omega)=0$.
Furthermore ${\bf P}^{12}_{\omega}$ contains a quadric
$Q_{p}$ for a point $p\in l$ if and only if $M_{l}(p)=0$.
\end{lemma}

\medskip

\proof Clear.
\qed

For the first part of the proposition we let $l$ be a $5$-secant line
to $Y$, and let $M(l)$ be the matrix of $\Phi$ restricted to $l$.
Assume first that $L$ is a non-isotropic $5$-secant line to $Y(3)$.  Then
$M(l)=(a_{1},\ldots,a_{5})$ has rank 2 in at 5 points.  But this is
possible only if $M(l)$ has rank 2 on all of $l$.  For a
general $\Phi$ this is not the case of any non-isotropic line.
If $l$ is an isotropic 5-secant line to $Y^s$ then again $M(l)$
must have rank 2 in 5 points on the line.   This is precisely the
case when $l$ is contained in $Y(3)$ or when $a(l)=0$.

For the last part of the proposition, consider first the set of
Lagrangian planes that pass through a point $p\in Y$.
In $X$, these planes are parameterized by a conic section, this is
\ref{prop:quadric bundle}.  In ${\bf P}(V)$, these planes form a pencil on
a $3$-dimensional quadric $Q$ in the hyperplane $H_{p}$ defined by
$L_{\alpha}(p)$.   The conic section $F_{p}$ in $Y$ through $p$
lies in a  Lagrangian plane whose point in $\Sigma$ is not in $X$.
In particular its plane is contained in $H_{p}$, but not in the
quadric $Q$.  The
intersection $Q\cap Y$ is therefore a curve of
degree at most $14$.  But every Lagrangian plane is $12$-secant,
unless the curve is rational it has degree at least $14$.  Therefore
we have equality.
When the conic in $X$ degenerates into two lines, the quadric
degenerates
into two $3$-spaces, which each intersects $Y$ in a curve of degree $7$
that obviously is a section of the conic bundle.
Let $C$ be such a curve on $Y$, let $F$ be conic, and let $H$ be a
hyperplane section.  Then $C$ has genus $3$, so by adjunction
$C^2+C\cdot K=C^2+C\cdot (-H+8F)=C^2-7+8=4$, i.e. $C^2=3$.  If
$D\subset Y$ is a section with $D^2\leq 2$, then $D\cdot K=D\cdot
(-H+8F)=-D\cdot H+8\geq 2$, i.e. $D\cdot H\leq 6$.  Since $Y$ is a
conic bundle the degree of the sections all have the same parity.  So
the degree of $D$ is $5$ or $3$. But $D$ is a smooth curve of genus
$3$, so this is impossible.
\qed

By \ref{fiber2}, \ref{prop:quadric bundle} it is clear that
the 3-fold $X$ parameterizes
$12$-secant planes to
$Y$, while $Y$ parameterizes conic sections contained
in $X$. 
\begin{proposition}\label{conics10}
      For the general smooth $3$-dimensional linear section $X$ of
$\Sigma$ the vertex variety
      $Y$ coincides with the Hilbert scheme $Hilb_{2l}(X)$ of
      conic sections on $X$.  \end{proposition}
      \proof  It remains to show that any conic section on $X$ is a
      plane section of the $3$-dimensional quadric $Q_{y}\subset\Sigma$
      of a point $y\in Y$.  But any conic section in $\Sigma$
      corresponds to one of the two pencils of planes on
      a $3$-dimensional quadric of
      rank $4$.  But these planes all have a common point $y$, so
      $Q_{y}\cap X$ is a conic, and $y\in Y$. \qed

\begin{proposition}\label{planes10}
     The general smooth $3$-fold linear section $X$ of $\Sigma$ coincides with 
     the set of planes in ${\bf P}(V)$ that intersect the vertex variety $Y$
     in a subscheme of length $12$.\end{proposition}
     \proof It remains to prove that a plane $P$ that intersects $Y$ 
     in a finite subscheme of length at least $12$ is Lagrangian and 
     belong to $X$.  
     Notice first that the $4$-space $\Pn 
     4_{y}$ intersects $<X>=\Pn {10}$ in a line when $y\in {\bf 
     P}(V)\setminus Y$, while the intersection is a plane when $y\in 
     Y$.  Since $Q_{y}=\Pn 4_{y}\cap \Sigma$ is a quadric threefold 
     that contains 
     lines but not planes, the intersection $Q_{y}\cap X$ may be of 
     three different kinds, namely 
     \par
     i) $Q_{y}\cap X$ is a scheme of length 
     $2$,   
    \par
     ii) $Q_{y}\cap 
     X$ is a line or 
     \par
     iii) $Q_{y}\cap X$ is a (possibly 
     singular) conic section. 
     
     \noindent Each plane in ${\bf P}(V)$ defines a hyperplane 
     section of $X$ consisting of those Lagrangian planes of $X$ 
     that intersect $P$.  Similarly to the points $y\in{\bf P}(V)$ are the following three possibilities 
     for a plane $P\subset {\bf P}(V)$:
     \par
     I)  There is a curve of points in $P$ through which passes 
     infinitely planes from $X$. 
     \par
     II)  Through only finitely many 
     points in $P$ pass infinitely many planes from $X$, but through 
     every  
     $y\in P$ there is a plane of $X$ that intersects $P$ along a line 
     through $y$.  
     \par
     III) Through only finitely many 
     points in $P$ pass infinitely many planes from $X$, and only 
     finitely many planes of $X$ intersect $P$ in a line. 
     \par\noindent
     Recall that $Y$ is a minimal conic bundle over a plane quartic curve, 
     whose sections of minimal selfintersection are curves $C$ of degree 
     $7$.  In particular the only plane curves on $Y$ are the smooth 
     conics.  Each minimal section $C$ lies in a $\Pn 3_{l}$ of a 
     $5$-secant line $l$ to $Y$.  In fact the scheme of intersection 
     $l\cap Y$ is the residual intersection to $C$ in $\Pn 3_{l}\cap Y$.    
     Assume that $P$ is a plane through a $5$-secant line $l$ that intersects $Y$ in a 
     scheme of length $12$.  If $P$ is not contained in $\Pn 3_{l}$, 
     i.e. it is not a plane of $X$, then the span $<P\cup \Pn 3_{l}>$ 
     is a hyperplane that intersects $Y$ in $C\cup C'$, where $C'$ is a section of 
     degree $9$.  Since $P\cap C'$ contains a scheme of length $12$, 
     we already have a contradiction.  
     Thus we may give a more precise description of when the different 
     cases above occur.\par\noindent
     For points case iii) occurs, of course, when $y\in Y$, while case ii) occurs when $y$ 
     lies outside $Y$ but on a $5$-secant line $l$ to $Y$.  For 
     planes case I) occur precisely when $P$ intersects $Y$ in a 
     conic or passes through a $5$-secant line $l$.  In the former 
     case $P$ is Lagrangian but does not belong to $X$.  The latter 
     case was argued above:  $P$ intersects $Y$ in a scheme of length
     $12$ only if $P$ is in $X$.
     For the case II) assume that infinitely many planes of $X$ meet $P$ in a 
     line. Since the planes of $X$ that intersects $P$ form a 
     hyperplane section of $X$, there still is a $2$-dimensional family 
     of planes of $X$ that meet $P$. Therefore there is a most one 
     plane of $X$ through a general point of $P$ that intersects $P$ in 
     a line, so the family of such lines of intersection form a line in 
     the dual space and  have some point $p$ in $Y$ in common.  This 
     point $p$ is obviously of type ii) or iii).  In case it is of 
     type ii), the plane $P$ and a $5$-secant line $l$ meet in $p$, 
     so they span a $\Pn 3$.  On the other hand this $\Pn 3$ 
     intersects $Y$ in a scheme of length at least $17$, so it 
     intersects $Y$ is a curve.  Any curve on $Y$ that span a $\Pn 3$ 
     is a minimal section, so $P$ must contain some $5$-secant line, 
     and belong to case I). In case $p$ is of type iii) 
    the planes of $X$ that meet $P$ in a line through $p$ form a conic 
    in $X$, and therefore forms 
one family of a quadric $3$-fold in ${\bf P}(V)$.  The plane $P$ is then a 
member of the other 
family, so $P$ and any plane of the conic intersects $Y$ in at most a 
$0$-dimensional $\Pn 3$ section of $Y$, i.e. in a scheme of length $16$.  
Since the plane of the conic meet $Y$ in a scheme of length $12$, and 
no line intersects $Y$ in a scheme of length $5$, the plane $P$ meets $Y$ 
in a scheme of length at most $9$, contrary to our assumption. 
The remaining  case III) follows from the following lemma. 
\begin{lemma}\label{lemmaplanes}
Let  $P\subset {\bf P}(V)$ be a plane that intersects $Y$ in a scheme of length $12$. 
Assume that only finitely many planes of $X$ meet $P$ in a line.  
Then $P$ belongs to $X$.\end{lemma}
\proof  Assume that $P$ does not belong to $X$.
If $P$ intersects a point of ii) we are done by the above 
argument. So the points of $P\cap Y$ are precisely the points of $P$ through 
which there pass infinitely many planes from $X$, i.e. a conic section 
of planes from $X$.  Since finitely many planes of $X$ 
meet $P$ in a line, there are (with multiplicity) through the general point in $P$ 
 precisely two planes of $X$, and they intersect $P$ 
precisely in the point.  Let $C\subset P$ be a curve and consider the 
curve $D$ of planes in $X$ that meet $C$. If $C$ intersects $Y$, 
there are conic sections in $X$ that are components of $D$. We 
define $D_{C}$ to be the complement of the conics in $D$.  Since the 
curve $D$ for a general line $C$ has degree $8$, the degree of $D$ for 
a curve $C$ of degree $d$ is $8d$.  The curve $D_{C}$ has degree 
$8d-2e$, where $e$ is the length of the subsceme $C\cap Y$. 
Clearly, for any curve $C$ in $P$, the curve $D_{C}$ has positive degree.  
Furthermore if the degree of $D_{C}$ is $2$, then $D_{C}$ is a conic 
or two disjoint lines. In the first case the planes of $D_{C}$ fill a 
quadric that intersects $P$ in a conic.  But since there are two 
planes in $D_{C}$ through a general point on $C$, this is impossible 
unless $C$ is a line. 
In the second case the plane $P$ meets two $5$-secant lines $l$, so 
this is ruled out above.  Therefore if a curve $C$ of degree $d\geq 2$ in $P$ passes 
through a subscheme of length $a$ of $Y$, then $8d-2e\geq 4,$ or 
$e\leq 4d-2$.\par\noindent     
Consider now the linear system of quartic curves 
in $P$ that pass through $P\cap Y$.   If $C$ is a line resp. conic or cubic 
curve in $P$ intersecting $Y$ in a scheme of length $a$, then $a\leq 3$ resp. 
$a\leq 6$ or $a\leq 10$.  Therefore the linear system of quartic 
curves through $P\cap Y$ can have no fixed component.
Furthermore the general meber of this system has at most one singular 
point (supported in $P\cap Y$), and this point is at worst an ordinary double 
point.   Let $C$ be a general such curve.    
Then $D_{C}$ is a curve of 
degree $8$ in $X$.  Since $P$ is not in $X$, there are two planes 
form $X$ distinct from $P$ through a general point of $C$,
so $D_{C}$ is a double covering of $C$. If $D_{C}$ is 
nonreduced, its reduction is a curve of genus at least $2$ and degree 
$4$, i.e. $D_{C}$ is a plane quartic curve.  But $X$ is cut out by 
quadrics so this is impossible.  If $D_{C}$ is reduced its genus is at 
least $3$, by Hurwitz' formula.  A curve of degree $8$ and genus at 
least $3$ spans at most a $\Pn 5$.  If the genus is $3$, the residual 
to a canonical divisor in a hyperplane section is a divisor of 
degree $4$ that spans a plane.  This divisor moves in a pencil, so 
$D_{C}$  has a pencil of $4$-secant planes. If the genus is bigger, 
the span of $D_{C}$ is smaller, so $D_{C}$ has even more $4$-secant 
planes.  But any plane that meet $X$ in a scheme of length at least 
$4$ intersects $X$ in a curve.  Therefore the span of $D_{C}$ 
intersects $X$ in a surface.  On the other hand the Picard group of 
$X$ is generated by a hyperplane section, so we get a contradiction.  
\qed

\begin{remark}\label{type2}  Mukai shows in \cite{M4} Theorem 9.1 that
$X$ is the Brill-Noether locus of Type II:  Let $F$ be a rank 2
vector bundle on a plane quartic curve of odd degree, such that any
section $C$ on the associated $\Pn 1$ bundle has self intersection
$C^2\geq 3$, then the moduli space $M$ of rank 2 vector bundles $E$, such
that ${\rm det}E-{\rm det}F=K$ and rank${\rm Hom}(F,E)\geq 3$, is a Fano
$3$-fold of index $1$ and genus $9$.  More precisely, $M$ is isomorphic
to a $3$-fold linear section $X$ of $\Sigma$.  The plane quartic curve is
the $Sp(3)$-dual curve $\check F(X)$.  The $\Pn 1$ bundle associated
to $F$ is nothing but the vertex surface $Y$.  It would be
interesting to see the identification of the bundles $E$ directly in
the above picture.
\end{remark}
In our setup the result of Mukai has the following corollary:
\begin{corollary}  Fix a smooth plane quartic curve $C$.
      There is a one-one correspondence between the
isomorphism classes of ruled surfaces $S$ over $C$ with minimal
selfintersection of a section equal to $3$ and the $3$-dimensional linear sections
$X\subset\Sigma$ such that the $Sp(3)$-dual section
$\check{F}(X)\cong C$.   The ruled surface $S$ is isomorphic to the
vertex variety $Y$ of $X$.
\end{corollary}

Iskovskikh \cite{I3} described the $3$-folds $X$, and in particular
their rationality via the double projection from a line:

\begin{theorem}\label{dproj}(Iskovskikh)
Let $ X $ be a smooth Fano $3$-fold of genus $g = 9$
of $rank \ Pic = 1$ and of index $1$, and let $L \subset X$
be a line. Then the double projection ${\pi}_{2L}$
from $L$, defined by the non-complete linear system
$|h-2L|$ on $X$, sends $X$ birationally to
$\Pn 3$.
Moreover, on $X$ there exists a unique cubic section
$M \in |O_X(3h-7L)|$ swept out by the 1-dimensional family
${\cal C}_L$ of conics $q \subset X$
which intersects the line $L$, and this family  of conics is contracted
to a smooth
curve ${\Gamma} = {\Gamma}^3_7 \subset {\bf P}^3$ of genus $3$
and of degree $7$ which lies on a unique cubic surface
$S_3 \subset {\bf P}^3$.

The inverse birational map
${\psi} = {\pi}_{2L}^{-1}: {\bf P}^3 \rightarrow X$
is given by the non-complete linear system
$|7H - 2 {\Gamma}|$ on ${\bf P}^3$.
Moreover, the cubic surface $S_3$ is swept out by the
$1$-dimensional family ${\cal C}_{\Gamma}$ of conics
$q \subset {\bf P}^3$ which are $7$-secant to ${\Gamma}$,
and this family of conics is contracted to the line $L \subset X$.

In addition, for any smooth curve
$\Gamma = {\Gamma}^3_7 \subset {\bf P}^3$
which lies in a unique cubic surface
$S_3 \subset {\bf P}^3$, the rational
map ${\psi}$ defined by the non-complete
linear system $|7H - 2{\Gamma}|$ on ${\bf P}^3$
defines a birational map from ${\bf P}^3$
to a smooth Fano $3$-fold $X = X_{16}$
of $rank \ Pic = 1$ and of index $1$,
and on this $X$ there exists a line $L$ such that
${\psi} = {\pi}_{2L}^{-1}$ where ${\pi}_{2L}$
is the double projection from $L$.
\end{theorem}

\centerline{The birationality between $X_{16}$
and the projective $3$-space $\Pn 3$}
\centerline{defined by double projection
                            from a line $L \subset X_{16}$}

\begin{picture}(100,40)
\put(70,34){\makebox(0,0){${\scriptstyle {\phi}\ =\ |h-2L|}$}}
\put(35,28){\makebox(0,0){{\footnotesize line} $L \subset X_{16}$}}
\put(60,30){\vector(1,0){20}}
\put(100,28){\makebox(0,0){${\bf P}^3 \supset {\Gamma}^3_{7}$}}
\put(80,26){\vector(-1,0){20}}
\put(70,23){\makebox(0,0)
   {${\scriptstyle {\psi}\ =\ |7H-2{\Gamma}|}$}}
\put(30,17){\makebox(0,0)
   {{\footnotesize extremal curves} ${\cal C}_L$}}
\put(110,17){\makebox(0,0)
   {{\footnotesize extremal curves} ${\cal C}_{\Gamma}$}}
\put(30,12){\makebox(0,0){{\footnotesize extremal divisor}}}
\put(30,7){\makebox(0,0){$M_l \equiv 3H - 7L$}}
\put(110,12){\makebox(0,0){{\footnotesize extremal divisor}}}
\put(110,7){\makebox(0,0){$N_{\Gamma} \equiv 3h - {\Gamma}$}}
\end{picture}

\paragraph{}\label{nb}
We may use \ref{dproj} to construct an isomorphism
$f: \check{F}(X) \rightarrow {\Gamma}^3_7$.

Consider a line $L\subset X$.  By \ref{lines10} the pencil of Lagrangian
planes in $\Pn 5$ corresponding to this line each intersects the
vertex surface $Y$
in a finite scheme of length $12$, and the common line is a $5$-secant.
The $\Pn 3$ union of these
planes cut $Y$ in a section $C_{L}$ of $Y(3)$ of degree  $7$.  On the
other hand $Y$ parameterizes conic sections on $X$.
The curve $C_{L}$ parameterizes a family of conic sections on $X$ which
intersects $L$.  Since $Y$ is a fiber bundle over the plane
quartic curve $\check{F}(X)$, we get an induced map from a family
of conics on $X$ that meet $L$ onto $\check{F}(X)$.  On the other hand
this family of conics is also fibered over $\Gamma$, so we get an
induced map from $\Gamma$ to $\check{F}(X)$.  Since they are both of
genus $3$, this must be an isomorphism.

\medskip

\begin{corollary}
The principally polarized intermediate jacobian $J_X$ of $X = X_{16}$
is isomorphic to the jacobian $J_{\check{F}(X)}$ of the plane quartic
curve $\check{F}(X)$, the $Sp(3)$-dual to
$X$.
\end{corollary}

\proof Indeed $J_X \cong J_{\Gamma}$
since ${\pi}_{2L}^{-1}: {\bf P}^3 \rightarrow X$
is a composition of a blow-up of $\Gamma$,
a flop and a blow-down of a divisor to $L \cong {\bf P}^1$
(see \cite{I3} and Corollary 9.7 of \cite{M4}).
Since $\Gamma \cong \check{F}(X)$ then $J_X \cong J_{\Gamma} \cong
J_{\check{F}(X)}$.
\qed

\section{Rank two vector bundles on linear sections}\label{sec:vector}

We turn to the main application of our study of $\Sigma$. On each
nodal hyperplane section of $\Sigma$ we will construct a rank two vector
bundle with $6$ global sections.

As before we use the following notation:
For a point $\omega\in \Pnd {13}$ we consider the hyperplane $\Pn
{12}_{\omega}$
and the hyperplane section $H_{\omega}=\Pn {12}_{\omega}\cap \Sigma$.  When
$\omega\in \check{F}\setminus\check{\Omega}\subset \Pnd {13}$ then
$u(\omega)\in \Sigma$
is the pivot of $\omega$ on $\Sigma$, and when $\omega\in
\check{\Omega}\setminus\check{\Sigma}\subset \Pnd {13}$,
then $Q_{\omega}\subset \Sigma$ is the smooth quadric surface of pivots of
$\omega$ on $\Sigma$ and also the singular locus of $H_{\omega}$.

\subsection{The projection of $H_\omega$ from the pivot
$u(\omega)$}\label{subsec:projection}

\paragraph{}\label{ha}
Let ${\bf P}^{12}_\omega \subset {\bf P}^{13}$ be the hyperplane
defined by $\omega \in \check{F} \subset \check{\bf P}^{13}$, let
${\pi}_u: {\bf P}^{12}_\omega \rto \overline{\bf P}^{11}_\omega$ be the
projection from $u = u(\omega) \in \Sigma$, and let
the variety $\overline{H}_\omega \subset \overline{\bf P}^{11}_\omega$ be the
proper ${\pi}_u$-image of the hyperplane section
$H_\omega = \Sigma \cap {\bf P}^{12}_\omega \subset {\bf P}^{12}_\omega$.

Let ${\sigma}: H'_\omega \rightarrow H_\omega$
be the blow-up of $u \in H_\omega$, and $\psi:
H'_\omega\to\overline{H}_\omega $ the projection into $\overline{\bf
P}^{11}_\omega$.
By \ref{singha} (i), $u = u(\omega)$ is an ordinary double point
of $H_\omega$, therefore the exceptional divisor
$Q' = {\sigma}^{-1}(u) \subset H'_\omega$ of ${\sigma}$
is isomorphic to a smooth $4$-dimensional quadric,
i.e. $Q' \cong G(2,{\bf C}^4)$.

\medskip

The projection ${\pi}_u$ contracts the tangent cone
   $C_u\subset {\bf P}^6_u= T_u{\Sigma}$ at $u$ to a Veronese surface $S_{u}$
(cf \ref{prop:tangsp}).

Since the exceptional divisor $Q' \subset H'_\omega$
is isomorphic to the projectivized tangent cone to $H_\omega$ at $u =
u(\omega)$,
the strict transform $C'_{u}$ of $C_{u}$ in $H'_\omega$ intersects
$Q'$ in a Veronese surface.
By \ref{singha}(i), \ $Q'$ is isomorphic to a smooth
$4$-dimensional quadric; therefore the isomorphic image
$\overline{Q} = {\psi}(Q') \subset \overline{H}_u$
is a smooth $4$-dimensional quadric containing the
surface $S_u$.

\vspace{0.5cm}

\begin{picture}(100,40)

\put(40,15){\makebox(0,0){$u \ \ \in$}}
\put(50,15){\makebox(0,0){$H_\omega$}}
\put(50,10){\makebox(0,0){$\cup$}}
\put(50, 5){\makebox(0,0){$C_u$}}
\put(55,15){
\vector(1,0){10}}
\put(80,15){\makebox(0,0)
             {$\overline{H}_\omega \supset \overline{Q} \supset S_u$}}
\put(40,35){\makebox(0,0){$Q' \ \subset \ $}}
\put(55,35){\makebox(0,0){$H'_\omega \supset C'_u$}}
\put(37,30){\vector(0,-1){10}}
\put(50,30){\vector(0,-1){10}}
\put(55,30){\vector(1,-1){10}}
\put(69,30){\vector(3,-2){13.5}}
\put(47,25){\makebox(0,0){${\sigma}$}}
\put(66,25){\makebox(0,0){${\psi}$}}
\put(60,17){\makebox(0,0){${\pi}_u$}}
\put(10,25){\makebox(0,0){{\rm The projection $\pi_{u}$}}}
\end{picture}

\vspace{0.5cm}

Since $\overline{H}_\omega$ is a birational
projection of $H_\omega$ from its double point $u = u(c)$, the degree
$deg \ \overline{H}_\omega = deg \ H_\omega - 2 = 14$.
Let $L$ be the hyperplane divisor on $H_\omega$.   We denote by $L$
also the pullback $
{\sigma}^*L$ and let $L'$ be the strict transform
of the general hyperplane divisor which passes
through the point $u$, i.e. $L'\equiv L-Q'$.
We summarize some further properties of the morphism $\psi$.
\medskip

\begin{lemma}\label{overlineH}

$\overline{H}_\omega \subset \overline{\bf P}^{11}_\omega$ is a
$5$-fold with singularities at most on the surface $S_u \subset \overline{Q}$.
$H'_\omega$ is a smooth $5$-fold, the morphism $\psi:H'_\omega\to
\overline{H}_\omega$ contracts the codimension 2 subvariety $C'_{u}$
to the surface $S_u \subset \overline{Q}$ and is an isomorphism
outside $C'_{u}$.  In particular, $\overline{H}_\omega$ has singularities
at most on the surface $S_u$.      The canonical divisor on
$H'_\omega$ is $K_{H'_\omega}=-
3L'$, where $L'$ is the pullback of a hyperplane divisor on
$\overline{H}_\omega$.
\end{lemma}
\proof  Since $\psi$ induces the projection from $u$, the divisor
$L'=L-Q'$ is the pullback of a hyperplane divisor on $\overline{H}_\omega$.
It remains only to compute the canonical divisor.
The canonical divisor on $\Sigma$ is $K_{\Sigma} \equiv -4H$, so by adjunction
the canonical divisor $K_{H_\omega} = -3L$.
Since
${\sigma}: H'_\omega \rightarrow H_\omega$ blows up the double point
$u \in H_\omega$,
the canonical divisor $$K_{H'_\omega}
\equiv {\sigma}^*K_{H_\omega} + (dim \ H_\omega \setminus
2)Q'\equiv{\sigma}^*(-3L) + 3Q' \equiv  - 3L'.$$

   \qed

\paragraph{}\label{hb}
We shall see in \ref{thm:hsection} below that the 5-fold
$\overline{H}_\omega$ is in fact
a linear section of the Grassmannian $G(2,6) \subset {\bf P}^{14}$
with a special codimension 3 subspace in ${\bf P}^{14}$.

Let $X^{6-k}={\bf P}^{13-k}\cap \Sigma$ be a smooth
$(6-k)$-dimensional linear section of
$\Sigma$ as in \ref{linear sections}. Let $\omega \in
\check{F}(X^{6-k})-\check{\Omega}(X^{6-k})$.  Then by \ref{singP} the pivot
$u(\omega)$ is not contained in ${\bf P}^{13-k}$.
Let ${\bf P}^{14-k}_{\omega}$ be the subspace of ${\bf P}^{13}={\bf 
P}(V(14))$ spanned
by ${\bf P}^{13-k}$ and $u(\omega)$, and let $\Pi^{k-2}_{\omega}=({\bf
P}^{14-k}_{\omega})^{\perp}$.  Clearly $\Pi^{k-2}_{\omega}$ is a
linear subspace of $\Pi^{k-1}$ of codimension one.
Denote by $W^{7-k}_{\omega}$ the intersection

$$W^{7-k}_{\omega} = \Sigma \cap {\bf P}^{14-k}_{\omega}.$$

Since the $(6-k)$-fold $X^{6-k}$ is a proper linear section of
$\Sigma$ and
   a linear section of $W^{7-k}_{\omega}$ with the codimension $1$ subspace
${\bf P}^{13-k} \subset {\bf P}^{14-k}_{\omega}$,
the dimension $dim \ W^{7-k}_{\omega} = 7-k$.
Furthermore, by \ref{singP}, the pivot $u({\omega})$ is a singular
point of $W^{7-k}_{\omega}$.

Consider now the projection ${\pi}_{u(\omega)}: {\bf P}^{12}_\omega
\rightarrow \overline{\bf
P}^{11}_\omega$
from the pivot $u(\omega)$.
Since $u(\omega) \not\in {\bf P}^{13-k} = <X^{6-k}>$,
the restriction of
${\pi}_{u(\omega)}$ to ${\bf P}^{13-k}$ is an projective-linear isomorphism
onto $
\overline{\bf P}^{13-k}_\omega := {\pi}_{u(b)}({\bf P}^{13-k})$,
in particular
${\pi}_{u(b)}: X^{6-k} \rightarrow
\overline{X}^{6-k}_\omega
= {\pi}_{u(b)}(X^{6-k}) \subset \overline{\bf P}^{13-k}_b$
is a projective-linear isomorphism.

Since ${\bf P}^{13-k} \subset {\bf P}^{14-k}_\omega$ is a hyperplane,
and since $u(\omega) \in {\bf P}^{14-k}_\omega$,
the projection ${\pi}_{u(\omega)}$ maps ${\bf P}^{14-k}_\omega$
onto $\overline{\bf P}^{13-k}_\omega$.
The pivot point $u(\omega)$ is a quadratic singularity of
$W^{7-k}_\omega$, therefore
the proper ${\pi}_{u(\omega)}$-image
$\overline{W}^{7-k}_\omega$
of $W^{7-k}_\omega$ will contain a quadric
$\overline{Q}^{6-k}_\omega \subset \overline{Q}$
of dimension $6-k$;
under the condition $2 \le k \le 5$.

We will show in \ref{thm:hsection} that the projection
${\pi}_{u(\omega)}$ sends the hyperplane section
$H_\omega = \Sigma \cap {\bf P}^{12}_\omega$
to a codimension $3$ linear section
$\overline{H}_\omega = \overline{\bf P}^{11} _\omega \cap G(2,{\bf C}^6)$
containing a smooth 4-fold quadric
$\overline{Q} = G(2,{\bf C}^4)$ for some ${\bf C}^4 \subset {\bf C}^6$.

Thus $\overline{X}^{6-k}_\omega$ is a subvariety of the linear
section $\overline{W}^{7-k}_\omega$ of $ G(2,{\bf C}^6)$, a linear
section that contains a $6-k$-dimensional quadric.

In the rest of this section these observations lie behind a description of
embeddings of linear sections of $\Sigma$ into $ G(2,{\bf C}^6)$, or
what amounts to the same, a description of a family of rank two
vector bundles on linear sections of $\Sigma$ with $6$ global sections.

\subsection{Del Pezzo and Segre threefolds}\label{DPS threefolds}

We define special Del Pezzo threefolds and identify them with
projections of Segre threefolds.  We shall show later that these are
subcanonical varieties on $\overline {H_\omega}$, and thus they are zero
loci of sections of a rank 2 vector bundle (via the Serre construction).

\paragraph{}\label{ia}
Let $\Gr 25 \subset {\bf P}^9 =
{\bf P}({\wedge}^2 \ {\bf C}^5)$ be the
Grassmannian of lines in ${\bf P}^4 = {\bf P}({\bf C}^5)$,
and let $\Gr 52\subset \check{\bf P}^9 =
{\bf P}({\wedge}^2 \ \check{\bf C}^5)$
be the Grassmannian of lines in the dual space
$\check{\bf P}^4 = {\bf P}(\check{\bf C}^5)$.
Any plane ${\Pi} \subset \check{\bf P}^9$
is the plane of linear equations of its zero-space
${\bf P}^6_{\Pi} \subset {\bf P}^9$.
The group $GL(5,{\bf C})$ acts by ${\wedge}^2$ on
the space $ \check{\bf P}^9$.
Therefore $GL(5,{\bf C})$ acts on
$\Gr 3{\wedge^2\check{\bf C}^5}$, the family of planes
in $\check{\bf P}^9$, and denote this action by
${\rho_{5}}$.

Let $U_0 \subset \Gr 3{\check{\bf C}^5}$
be the open set of these planes
$\Pi \subset \check{\bf P}^9$ such that
$\Pi \cap \Gr 52= \emptyset$.
The linear section $V_{\Pi} = {\bf P}^6_{\Pi} \cap \Gr 25$ is singular
if and only if it is contained in a tangent hyperplane since the
contact locus of a tangent hyperplane is a plane.  Therefore
$\Pi \in U_0$
\ {\it iff} \
$V_{\Pi} = {\bf P}^6_{\Pi} \cap \Gr 25$ is a smooth
Fano threefold of degree 5 and of index 2,
(i.e. $K_{V_{\Pi}} \equiv {\cal O}_{V_{\Pi}}(-2)$).
The action ${\rho_{5}}$ is transitive on $U_0$,
i.e. $U_0$ is an orbit of ${\rho_{5}}$ (Proposition 14
\S 5 in \cite{SK}).
Therefore all the $V_{\Pi}, {\Pi} \in U_0$ are conjugate to each other
by the action ${\wedge}^2$ of $GL(5,{\bf C})$ on ${\bf P}^9$
(see also \cite{I1} \S (6.5)).

\paragraph{}\label{ib}

Let $U_{xxx} \subset \Gr 3{{\wedge}^2 \ \check{\bf C}^5}$
be the subset of these planes
$\Pi \subset \check{\bf P}^9$ such that
$\Pi$ intersects $\Gr 52$ transversally
in exactly three points.
These three points of intersection cannot be collinear, since $\Gr 52$
is an intersection of quadrics.

As above, $U_{xxx}$ is an orbit of $\rho_{5}$,
and all the $V_{\Pi}, {\Pi} \in U_{xxx}$ are conjugate to each other
by the action ${\wedge}^2$ of $GL(5,{\bf C})$ on ${\bf P}^9$.

We call the unique threefold $V_{\Pi}$, ${\Pi} \in U_{xxx}$
the {\it Del Pezzo threefold of type $xxx$}.


\bigskip


Next, we turn to Segre threefolds.  For this we first make a slight
detour to $2-$forms on even dimensional spaces
and prove:

\begin{proposition}\label{diag} \
Let $V = {\bf C}^{2n}$ and let $\alpha, \alpha^{\prime}\in \wedge^2V$ be two
general $2-$vectors.  Then there exists a unique (up to scalars) $n-$ tuple
$\gamma_{1},\ldots, \gamma_{n}$ of
$2-$vectors of rank $2$ such that both $\alpha$ and $\alpha^{\prime}$ are
linear combinations of the $\gamma_{i}$.
\end{proposition}

\begin{remark}\label{n-secant}
      The proposition may be reformulated in terms of multisecant spaces to
the Grassmannian $G(2,2n)$ of lines in ${\bf P}^{2n-1}$ embedded in
Pl\"ucker space: A general line in the Pl\"ucker space is contained in a
unique $n$-secant $(n-1)$-space to $G(2,2n)$.
\end{remark}
\medskip
{\sl Proof of \ref{diag}.}  First we prove uniqueness.  Since the pair
$\alpha, \alpha^{\prime}$
of $2-$vectors is general we may suppose that they both have rank $2n$
and that
$$\alpha=\sum_{i=1}^n\gamma_{i},\quad{\rm and}\quad
\alpha^{\prime}=\sum_{i=1}^n\lambda_{i}\gamma_{i},$$
where the $\lambda_{i}$ are pairwise distinct coefficients.  Let
$$\beta_{i}=\lambda_{i}\alpha-\alpha^{\prime},\quad i=1,\ldots , n.$$
Then the $\beta_{i}$ are the precisely the $2-$vectors of the pencil
generated by $\alpha$ and $\alpha^{\prime}$ that have rank less than
$2n$.  Furthermore their rank is exactly $2n-2$
since $\lambda_i \not= \lambda_j$ for $i \not= j$.
Therefore each $\beta_{i}\in
\wedge^2V_{i}$ for a unique rank $2n-2$ subspace $V_{i}\subset V$.
Let $U_{j}=\cap_{i\not= j}V_{i}$.  Then $U_{j}$ is $2$-dimensional and
$\gamma_{j}$ is a nonzero $2-$vector that generates the subspace
$\wedge^2U_{j}\in\wedge^2V$,
so the $2-$vectors $\gamma_{i}$ are determined uniquely by the pencil
generated by $\alpha$ and $\alpha^{\prime}$.

For existence we use a dimension argument.  On the one hand,
the Grassmannian of lines in the Pl\"ucker space ${\bf P}(\wedge^2V)$
has dimension $2(n(2n-1)-2)=4n^2-2n-4$.  On the other hand the
Grassmannian $G(2,2n)$ has dimension $4n-4$, so the family of lines
contained in $n-$secant $(n-1)$-spaces to $G(2,2n)$ has dimension at
most $(4n-4)n+2(n-2)=4n^2-2n-4$.  By the uniqueness argument, a general
such line lies in a unique $n-$secant $(n-1)$-space, so the two
dimensions actually coincide.  Since there is an obvious inclusion of
the latter into the former, and the former is irreducible, we may
conclude.\qed

This result leads to a simple description of Segre $n-$folds
$$X_{n}=\Pn 1\times \Pn 1\times\cdots\times\Pn 1$$
in the
Grassmannian $G(n,2n)$.

\begin{proposition}\label{segre nfold}
Let $V= {\bf C}^{2n}$ and let $\alpha, \alpha^{\prime}\in
\wedge^2V^\ast$ be two
general $2-$forms on $V$.  Then the set of common Lagrangian $n-$spaces of
$V$ with respect to the forms $\alpha$ and $\alpha^{\prime}$ form a
Segre $n-$fold ${\bf
P}^{1}\times {\bf P}^{1}\times\cdots\times{\bf P}^{1}$ in the
Grassmannian $G(n,2n)$.
Let
$\gamma_{1},\ldots, \gamma_{n}$ be the unique $n-$tuple of
$2-$forms of rank $2$ such that both $\alpha$ and $\alpha^{\prime}$ are
linear combinations of the $\gamma_{i}$, and let $U_{i}\subset V$
be the $2n-2$-dimensional
kernel of $\gamma_{i}$.  Then the common Lagrangian $n-$spaces $U$ with
respect to the forms $\alpha$ and $\alpha^{\prime}$ are precisely
the $n-$spaces that intersects each $U_{i}$ in
a $n-1$-space. Equivalently, if $W_{i}=\cap_{i\not=j}U_{j}$, then
$W_{i}$ is 2-dimensional, and a $n-$space is Lagrangian with respect to
$\alpha$ and $\alpha^{\prime}$ if and only if it has a nontrivial
intersection with
each $W_{i}$.
\end{proposition}

\proof First, the two characterizations of Lagrangian $n-$spaces are
clearly equivalent.  Furthermore, the last one describe a family of
$n-$spaces that clearly form a Segre $n$-fold, since every $n$-space
intersects each line ${\bf P}(W_{i})$ in a unique point.

A $n-$space $U$ that intersects each $U_{i}$ in
a $n-1$-space is clearly isotropic with respect to each $2$-form
$\gamma_{i}$, therefore also Lagrangian for $\alpha$ and $\alpha^{\prime}$.

On the other hand it is a straightforward exercise in Schubert calculus
to show that the set of common Lagrangian $n-$spaces for $\alpha$ and
$\alpha^{\prime}$ is $n-$dimensional of degree $n!$ i.e. the degree of
the Segre
$n-$fold.   Since we have an inclusion the conclusion follows.
\qed

We turn back to the case $n=3$:

\begin{lemma}\label{xxx}
Let $X = X_3 \subset {\bf P}^7$ be the Segre threefold,
and let $u \in X$. Then the projection
$\overline{V}$ of $X$ from $u$ is a Del Pezzo threefold
of type $xxx$.

Conversely, the Del Pezzo threefold $\overline{V}$ of type $xxx$
is a projection of the Segre threefold $X$ from a point
$u \in X$.
\end{lemma}

\proof Consider the blowup $X'\to X$ centered at $u$, and let $E$ denote the
exceptional divisor.
Let $L_{i}$ be the pullback to $X$ of ${\cal O}_{{\bf P}^1}(1)$ on each
factor of $X$, and let $F_{X}=L_{1}\oplus L_{2}\oplus L_{3}$.  Let $s_{u}$ be
the unique global section of $F_{X}$ whose zero locus is $u$, and
let $F_{X'}$ be the pullback of $F_{X}$
$X'$. The pullback of $s_{u}$ to $F_{X'}$ vanishes on $E$, and
correspond to the unique nonvanishing section $s'$ of $F_{X'}(-E)$.
The exterior multiplication with $s'$ defines a surjective map
$\wedge^2F_{X'}(-E)\to \wedge^3F_{X'}(-2E)$ that fits into an exact sequence

$$0\to F_{0}\to \wedge^2F_{X'}(-E)\to \wedge^3F_{X'}(-2E)\to 0,$$
where $F_{0}$ is a rank 2 vector bundle on $X'$.  Notice
$\wedge^3F_{X'}(-2E)$ is the line bundle ${\cal O}_{X'}(H_{X}-2E)$
where $H_{X}$
is the pullback to $X'$ of the
hyperplane class on $X$ in the Segre embedding.  Furthermore
$c_{1}(\wedge^2F_{X'}(-E))=2H_{X}-3E$ so $c_{1}(F_{0})=H_{X}-E$. Similarly
one computes $c_{2}(F_{0})=(H_{X}-E)^2$  On the
other hand, $h^0({\cal O}_{X'}(H_{X}-2E))=4$, while s there
$h^0(\wedge^2F_{X'}(-E))=9$ so $h^0(F_{0}\geq 5$.  Therefore the
morphism defined by $H_{X}-E$, i.e. the projection of $X$ from $u$, maps
$X'$ into $\Gr 25$.  It is a threefold of degree 5 that spans a $\Pn
6$, so it is a linear section of $\Gr 25$.

$X$ contains three quadric surfaces that meet pairwise
along a line through $u$.  Thus the projection $\overline{X}$, i.e.
the image of $X'$ contains three
planes that meet pairwise in three points.  Evidently  these points are
precisely the singularities of $\overline{X}$.  Thus $\overline{X}$
is a Del Pezzo threefold of type $xxx$.

The converse is clear by the transitivity of $\rho_{5}$.  \qed

\begin{lemma}\label{segre threefolds}
Let $X = X_3 \subset \Sigma$ be a Segre threefold,
and let $H$ be a hyperplane section of $\Sigma$ that contains $X$. Then
$H$  is singular, i.e. a tangent hyperplane section
to $\Sigma$, at some point of $X$.
\end{lemma}

\proof For each point $u$ on $X$ there is a ${\bf P}^2$ of hyperplanes
tangent to $\Sigma$ at $u$ that contain $X$.  Since the general such
hyperplane is tangent to $\Sigma$ at $u$ only, there is altogether a
$5$ dimensional family of tangent hyperplanes that contain $X$.  But
$X$ spans a ${\bf P}^7$ so there is a ${\bf P}^5$ of
hyperplanes that contain $X$.  Therefore the two sets must coincide
and the lemma follows.

\qed

\subsection{A rank 2 vector bundle on singular
hyperplane sections}\label{subsec:vb}

\paragraph{}\label{ja}
Recall the universal sequence of vector bundles on $\Sigma$, the
restriction of the universal sequence on $G=\Gr 3V$:

\medskip

\centerline{
$0\to U\to V\otimes {\cal O}_{\Sigma}\to Q\to 0$,
}

\medskip

\noindent
where $U$ is the universal subbundle, and $U^{*}\cong Q$, by the
natural map induced by $\alpha$.
Any global section of the rank 3 bundle $\wedge^2 U^{*}$ comes from a
2-form $\alpha^{\prime}\in \wedge^2V^{*}$.
In the previous section we saw that if the 2-form $\alpha^{\prime}$ is
general, then the zero-locus $X=Z(\alpha^{\prime})$ is a Segre
3-fold.  In fact the characterization \ref{segre nfold} yields a
straightforward argument that any Segre 3-fold in $\Sigma$ is the zero-locus
a section of $\wedge^2 U^{*}$.  Since it is not essential we leave
this out here.

  From \ref{segre threefolds} we know that any hyperplane that contains
the Segre 3-fold
$X$ is tangent to $\Sigma$.  So we fix a hyperplane $\Pn {12}_{\omega}$ that
contains $X$,
and assume that it is tangent only at $u\in X$. Then $X$ has codimension 2
in this hyperplane section, and still it is the zero-locus of a $2$-form
restricted to the hyperplane.  Consider the blowup $H'_{\omega}$ of
the hyperplane
section $H_{\omega}=\Sigma\cap \Pn {12}_{\omega}$ in the singular point $u$.
Let $Q'_{\omega}$ be the exceptional divisor on $H'_{\omega}$.  It
is isomorphic to a 3-dimensional quadric (when $\omega$ is general).

In the notation of \ref{overlineH} the canonical bundle $K_{X}={\cal
O}_{X}(-2L)$.  On the strict transform $X_{\omega}$ of $X$, the
canonical bundle is therefore $$K_{X_{\omega}}={\cal
O}_{X_{\omega}}(-2L+2Q'_{\omega})={\cal
O}_{X_{\omega}}(-2L'),$$ so $X_{\omega}$ is subcanonical with respect
to the hyperplane line bundle $L'$ induced by the projection from
$u$. Therefore,
by the Serre construction, $X_{\omega}$ is the zero-locus of a rank $2$
vector bundle on $H'_{\omega}$.
The aim of this section is to identify this vector bundle.  To
construct it one may apply the Serre construction starting with $X_{\omega}$.
We will use a different, and for our purposes, more direct argument,
similar to the one used in the proof of \ref{xxx}.

Let $$\wedge^2 U^{*}_{\omega}=\wedge^2 U^{*}\otimes {\cal
O}_{H'_{\omega}},$$ and consider the twisted bundle $$\wedge^2
U^{*}_{\omega}(-Q'_{\omega}).$$ Notice that
the section of  $\wedge^2 U^{*}_{\omega}$,
given by the restriction and pullback of the
section $\alpha^{\prime}$  vanishes on $Q'_{\omega}$.  Therefore it
corresponds to a section $\alpha^{\prime}_{\omega}$ of $\wedge^2
U^{*}_{\omega}(-Q'_{\omega})$.  The vector bundle   $\wedge^2
U^{*}_{\omega}(-Q'_{\omega})$ has rank
3, while the zero-locus $X_{\omega}$ of the section
$\alpha^{\prime}_{\omega}$ has
codimension 2.  We will show that $\alpha^{\prime}_{\omega}$ is a
section of a rank 2 subbundle of $\wedge^2 U^{*}_{\omega}(-Q'_{\omega})$.
To do this we consider maps of
vector bundles:
$$\wedge^2 U^{*}_{\omega}(-Q'_{\omega})\to {\cal
O}_{H'_{\omega}}(L-2Q'),$$
where $L$ is the pullback
of the general hyperplane divisor on $H_{\omega}$.
Let $U\subset V$ be the Lagrangian $3-$space represented by $u\in
\Sigma$ and let $U^{\bot}=L_{\alpha}(U)\subset V^{*}$.   Then
any element $x\in U^{\bot}$ induces by exterior multiplication such a
map:
$$m_{x}:\wedge^2 U^{*}_{\omega}(-Q'_{\omega})\to {\cal
O}_{H'_{\omega}}(L-2Q'_{\omega}).$$
The kernel of this map, we denote it by $E'_{x}$, is of course a torsion
free sheaf.  If the map $m_{x}$ is surjective, the kernel is even a
vector bundle of rank $2$.  Therefore $E'_{x}$ is our candidate for a
rank 2 vector bundle.
If we look at the stalks, we see that the multiplication by
$x$ is surjective outside the zero locus of $x$, i.e. outside the strict
transform $C(x)$ on $H'_{\omega}$ of the
quadric cone $Q_{x}\cap H_{\omega}$ with vertex at $u$.  Since
$Q_{x}$ is $3$-dimensional, $C(x)$ is a (rational) surface scroll
whose image in ${\overline H}_{\omega}$ is a conic section on the
Veronese surface $S_{u}$.     Thus we have an exact
sequence

$$0\to E'_{x}\to \wedge^2 U^{*}_{\omega}(-Q'_{\omega})\to {\cal
O}_{H'_{\omega}}(L-2Q'_{\omega})\to {\cal
O}_{C(x)}(L-2Q'_{\omega})\to 0.$$

Outside $C(x)$ the kernel sheaf $E'_{x}$ is a rank 2 vector bundle.
    This
will be enough for our purposes at this point, but eventually we will
show that $E'_{x}$ is a subscheaf of a bundle $E_{x}$ that coincides
with $E'_{x}$ outside $C_{x}$.

The problem is to get $h^0(E'_{x})=6$.  By restriction and pullback
from $\Sigma$ there is a natural surjection of sections
$$r_{\omega}:V^{*}(14)\cap (\wedge^2 U^{\bot}\otimes
U_{1}^{\bot}\oplus \wedge^3
U^{\bot})= Sym^2U^{*}\oplus \wedge^3
U^{\bot}\to H^0({\cal
O}_{H'_{\omega}}(L-2Q'_{\omega})).$$  The kernel of this map is
generated by $\omega\in Sym^2U^{*}$, so
$h^0({\cal
O}_{H'_{\omega}}(L-2Q'_{\omega})=6$.

Similarly there is a natural surjection $$U^{\bot}\otimes U_{1}^{\bot}\oplus
\wedge^2
U^{\bot}\to H^0(\wedge^2
U^{*}_{\omega}(-Q'_{\omega})).$$  Here the kernel is generated by
$\alpha$, while $U^{\bot}\otimes U_{1}^{\bot}\oplus
\wedge^2
U^{\bot}$ is $12$-dimensional, so $h^0(\wedge^2
U^{*}_{\omega}(-Q'_{\omega}))=11$.  Thus $h^0(E'_{x})=6$ only if the
map $m_{x}$ is not surjective on global sections.

We consider more carefully the image of
the map $m_{x}$ on global sections.  Notice that $r_{\omega}(\eta)$ for a form
$$\eta\in V^{*}(14)\cap (\wedge^2 U^{\bot}\otimes U_{1}^{\bot}\oplus
\wedge^3 U^{\bot})$$
is in the image of $m_{x}$ if and only if there exists a
2-form $\beta\in U^{\bot}\otimes U_{1}^{\bot}\oplus
\wedge^2 U^{\bot}$  and a 1-form $y\in U^{\bot} $ such that
$$\eta=\alpha\wedge y + \beta\wedge x.$$
The subspace of 3-forms of this kind in $\wedge^2 U^{\bot}\otimes
U_{1}^{\bot}\oplus
\wedge^3 U^{\bot}$ has dimension $9$, i.e. codimension $1$:  The
3-forms of the kind $\alpha\wedge y$ form a $3$-dimensional space,
while the $3$-forms of the kind $\beta\wedge x$, where $\beta$ varies,
form a subspace of dimension $7$. Since these two subspaces intersect
in $<\alpha\wedge x>$, the dimension of their sum is $9$.  The
intersection with $V^{*}(14)$ has codimension $3$, it is the
symmetrizer relations in $V^{*}(14)\cap \wedge^2 U^{\bot}\otimes
U_{1}^{\bot}$, so the subspace
$$U_{x}=\{\eta=\alpha\wedge y + \beta\wedge x|\eta\wedge
\alpha=0\}\subset V^{*}(14)\cap (\wedge^2 U^{\bot}\otimes U_{1}^{\bot}\oplus
\wedge^3 U^{\bot})$$
has dimension
6. The image of the map $m_{x}$ on global sections is nothing but the
projection of $U_{x}$ from the form $\omega$.  Therefore we have shown

\begin{lemma}\label{nsurjective}
     The exterior multiplication $$m_{x}:\wedge^2
U^{*}_{\omega}(-Q'_{\omega})\to {\cal
O}_{H'_{\omega}}(L-2Q'_{\omega})$$ is not surjective on global
sections if and only if $\omega$
is an element of $U_{x}$.
    \end{lemma}

\bigskip
Let $p=<v>\in {\bf P}(U)$ and $x=L_{\alpha}(v)\in U^{\bot}$. Let
$q_{\omega}$ be the quadratic form defined by $\omega$ on $U$ (cf.
\ref{subsec:rep}).
\begin{lemma}\label{nonsurjective}
     Let $\omega\in V^{*}(14)\cap (\wedge^2 U^{\bot}\otimes
U_{1}^{\bot}\oplus \wedge^3
U^{\bot})$.  Then $\omega\in U_{x}$ if and only if $q_{\omega}(v)=0$.
    \end{lemma}
\proof  First we assume that
$$\omega=\alpha\wedge y+\beta\wedge x.$$
The common zero-locus of the $2$-forms $\alpha$ and $\beta$ is then
contained in  $H_{\omega}$.
Therefore we may choose coordinates $(e_{i},x_{i}$ on $V$ such that
   $U=U_{0}=<e_{1},e_{2},e_{3}>$,  $U_{1}=<e_{4},e_{5},e_{6}>$
    and assume that
$$\beta=sx_{14}+tx_{25}+ux_{36}.$$
Since $\alpha\wedge \omega =(x_{14}+x_{25}+x_{36})\wedge\omega=0$,
then
$$\omega = b(x_{145}+x_{356})-c(x_{416}+x_{256})-a(x_{452}+x_{436}),$$
for suitable scalar coefficients $a,b,c$.

The quadratic form $q_{\omega}$ on $U$ is then (cf. \ref{subsec:rep})
$$q_{\omega}=bx_{1}x_{3}-cx_{1}x_{2}-ax_{2}x_{3}.$$
In the expression
$$\omega=\beta\wedge x+\alpha\wedge y,$$
   it is clear that
$$x,y\in
<x_{4},x_{5},x_{6}>=L_{\alpha}(U).$$
\noindent
Thus we may write $x=\beta_{4}x_{4}+\beta_{5}x_{5}+\beta_{6}x_{6}$ and
$y=\alpha_{4}x_{4}+\alpha_{5}x_{5}+\alpha_{6}x_{6}$.
A straightforward calculation gives the following solutions:
$$\alpha_{4}=a{{u+t}\over{u-t}}, \alpha_{5}=b{{u+s}\over{u-s}},
\alpha_{6}=c{{t+s}\over{t-s}}$$
and
$$\beta_{4}={a\over{u-t}},
\beta_{5}={b\over{u-s}}, \beta_{6}={c\over{t-s}},$$
thus $v={a\over{u-t}}e_{1}+{b\over{u-s}}e_{2}+{c\over{t-s}}e_{3}$ and
$q_{\omega}(v)=0$.

Conversely assume $q_{\omega}(v)=0$. Let $X$ be a Segre $3$-fold through
$u$ contained in $H_{\omega}$. Then we may assume that $X$ is the zero
locus of a $2$-form $\beta$.
Coordinates may therefore be chosen as above, and  $q_{\omega}(v)=0$
implies that
$$v={a\over{u-t}}e_{1}+{b\over{u-s}}e_{2}+{c\over{t-s}}e_{3}.$$
With
$$y=a{{u+t}\over{u-t}}x_{4}+b{{u+s}\over{u-s}}x_{5}+c{{t+s}\over{t-s}}
x_{6}$$
we get
$$\omega=\alpha\wedge y+\beta\wedge x.$$
   \qed

\begin{corollary}\label{6 sections}
Let $v\in U$, and $x=L_{\alpha}(v)\in U^{\bot}$.  Let $E'_{x}$ on
$H'_{\omega}$ be the
kernel sheaf of the map $m_{x}$ above.  Then $h^0(E'_{x})=6$
   if and only if $q_{\omega}(v)=0$, where
$H'_{\omega}$ is tangent at $u\in \Sigma$ and $q_{\omega}$ is the
quadratic form defined by $\omega$ on $U$.
\end{corollary}
\proof  Since $\wedge^2 U^{*}_{\omega}(-Q'_{\omega})$ has 11 sections,
$h^0(E'_{x})=6$ if and only if
$m_{x}$ is not surjective on global sections, i.e.
$\omega=\alpha\wedge y+\beta\wedge x$ for some $1$-form $y$ and
$2$-form $\beta$.  So the corollary follows from \ref{nonsurjective}.\qed

For each $v\in  U$ such that $\{q_{\omega}(v)=0\}$ we have contructed
a sheaf $E'_{x}$, where $x=L_{\alpha}(v)$, locally free of rank 2 outside
$C(x)$  on $H'_{\omega}$  with $h^0(E'_{x})=6$.
Each of them gives rise to a rational map of $H'_{\omega}$
into the Grassmannian $\Gr 26$.   This map is defined by sections of the
determinant line bundle of $E'_{x}$, whose first Chern class is given by:
$$c_{1}(E_{x})=c_{1}(\wedge^2 U^{*}_{\omega}(-Q'_{\omega}))-c_{1}({\cal
O}_{\Sigma_{\omega}}(L-2Q'_{\omega}))$$
$$=(2L-3Q'_{\omega})-(L-2Q'_{\omega})=L-Q'_{\omega}=L'.$$
If the natural map
$$\wedge^2H^0(E'_{x})\to H^0{\cal
O}_{H'_{\omega}}(L')$$
is surjective, then the map to $\Gr 26$ is nothing but the
projection of $H'_{\omega}$ from its singular point $u=u(\omega)$.
If it is not surjective
the image in $\Gr 26$ spans at most a $\Pn {10}$.  Since $H'_{\omega}$
is not a cone, the image is $5$-dimensional, and the intersection of
its span with
$\Gr 26$ is not proper.  It is now a straightforward but tedious to
check that no $\Pn {10}$ section of $\Gr 26$ can contain the image
of $H'_{\omega}$.  Therefore the image in $\Gr 26$ is precisely the projection
$\overline {H}_{\omega}$.  Furthermore
this map is independent of $x$.

\begin{theorem}\label{thm:hsection}
The projection of $H_{\omega}$ from the singular point is a
linear section of the Grassmannian $\Gr 26$.
\end{theorem}

\proof  What remains is to show that the image
$\overline{H}_{\omega}$ of $H'_{\omega}$
under the
projection is a linear section of $\Gr 26$.  Since $\Sigma$ has
degree $16$ and sectional genus $9$, the projection of
$H_{\omega}$ must have degree $14$ and sectional genus $8$.  It
is $5$-dimensional, spans a $\Pn {11}$ and is contained in $\Gr 26$.
Furthermore, it contains a $4$-dimensional quadric, the image of
the exceptional divisor on  $H'_{\omega}$.  Thus $\overline{H}_{\omega}$ 
has the same degree, sectional genus and codimension as $\Gr 26$.  So if the intersection
$$\Gr 26\cap <\overline{H}_{\omega}>$$

is different from $\overline{H}_{\omega}$, then this intersection is
6-dimensional and not proper.   But the only codimension 2 varieties in
$\Gr 26$ that are contained in a $\Pn {11}$ are those representing
special Schubert
cycles of codimension $2$. One is represented by the subvariety of
rank 2 subspaces that
intersect a given rank $3$ subspace and the other is represented by
Grassmannians $\Gr 25$.  The latter does not span a $\Pn {11}$,
while the former do.  In the former case
$\overline{H}_{\omega}$ is contained in a variety that is a
$\Pn 4$ scroll parameterized by a ${\bf P}^2$. Each $\Pn 4$ must
intersect $\overline{H}_{\omega}$ in a threefold.  But
$H_{\omega}$ is cut out by quadrics and contains only a
$1$-parameter family of threefold hypersurfaces (in fact the  linear
$3$-spaces that appear as projections of $Q_{p}\subset H'_{\omega}$ of
\ref{incidence}), so we have a contradiction and the theorem follows.
\qed

\medskip
The restriction and pullback to $H'_{\omega}$
of the universal rank $2$ quotient bundle on $\Gr 26$ is clearly a rank 2
vector bundle.  We denote it by $E_{\omega}$.
We immediately get:

\begin{corollary}\label{taut}
The sheaf $E'_{x}$ is a subsheaf of the restriction and pullback
$E_{\omega}$ to
$H'_{\omega}$ of the universal rank $2$ quotient bundle on $\Gr 26$.
The bundle $E_{\omega}$ is independent of
$x$, has $h^0(E_{\omega})=6$,
det$E_{\omega}={\cal O}_{H'_{\omega}}(L')$ and the zero scheme of a
general section is isomorphic to the strict transform of a Segre threefold
that passes through the singular point.
\end{corollary}
\proof Outside $C(x)$ the two sheaves $E'_{x}$ and
$E_{\omega}$ coincides. The zeros of a general section is precisely
the zeros of a $2$-form on $\Sigma$, i.e. the strict transform on
$H'_{\omega}$ of a Segre threefold that passes through the singular 
point $u$. Since
$C(x)$ has codimension $3$ the
corollary follows.\qed

Clearly $\overline{H}_{\omega}$ is a special linear section of
$\Gr 26$ since it contains a $4$-dimensional quadric, but a natural
question arises: Is a general $\Pn {11}$ section  of
$\Gr 26$ that contains a $4$-dimensional quadric the projection of a
singular section of the Lagrangian Grassmannian $\Sigma$.

We now prove that this is the case, and give another
characterization of these linear sections of $\Gr 26$.
We set $$Z = \overline{H}_\omega,$$
and notice that the projection of $H'_{\omega}$ is an isomorphism
on the exceptional quadric and outside the tangent cone at $u$.  Thus
it is singular at most along the image of the tangent cone, i.e. a
Veronese surface (inside the $4$-dimensional quadric).

First, since $\check{F} \setminus \check{\Omega}$ is an orbit for
$\rho$, all the
$\overline{H}_\omega$ are projectively equivalent
to the same 5-fold $Z$.

To fix notation, let $V\cong {\bf C}^6$ and let
   ${\bf P}^{14} =
{\bf P}({\wedge}^2 V)$.
Let $\Gr 2V$ be the Grassmannian
of 2-dimensional subspaces $U \subset V$,
and let

$$\Gr 2V \rightarrow {\bf P}^{14},\quad
U \mapsto {\bf P}({\wedge}^2 U)$$
be the Pl\"ucker embedding.

Let $\check{\bf P}^{14} = {\bf P}({\wedge}^2 {V^{*}})$ be the dual
space to $\Pn {14}$.  The space
${\wedge}^2 V^{*}$
is isomorphic to the space $Alt(V,V^{*})$
of skew-symmetric linear maps
$A: V \rightarrow V^{*}.$
Recall that the rank of $A$ is {\it even}.  The rank stratification is
given by the inclusions

$$\Gr V2= \Gr 2{V^{*}} \subset \check{Pf}
\subset \check{\bf P}^{14}$$
that is by the  Grassmannian variety parameterizing ${\bf C}^*$-classes
of maps $A \not= 0$
s.t. $rank(A) = 2$, and the Pfaffian cubic hypersurface parameterizing
${\bf C}^*$-classes
of maps $A$ s.t. $rank(A) \le 4$.

Let ${\Pi}^2 \subset \check{\bf P}^{14}$ be the plane
of linear equations that define $\Pn {11} \subset {\bf P}^{14}$,
such that the 5-fold $Z'=\Pn {11}\cap \Gr 2V$ contains a smooth
4-dimensional quadric $Q$.

Obviously, $Q = \Gr 2W \subset \Gr 2V$
for some $4-$dimensional subspace $W\subset V$.

\begin{lemma}\label{pfplane}
${\Pi}^2 \subset \check{Pf}$.
\end{lemma}

\proof The forms in $V^{*}$ that vanish on the $4-$dimensional
subspace $W\subset V$ is a rank $2$
subspace $W^{\perp}\subset V^{*}$.
Then
$$({\wedge}^2 \ W)^{\perp}
= V^{*} \wedge W^{\perp}\subset {\wedge}^2 V^{*}.$$
Therefore any $A \in ({\wedge}^2 \ W)^{\perp}$
is of rank at most $4$, i.e.
$${\bf P}(({\wedge}^2W)^{\perp}) \subset \check{Pf}.$$
Since $Q=G(2,W)\subset Z'$, the lemma follows.
   \qed

Fix $W\subset V$ a $4-$dimensional subspace, and let $\Pn 8_{W}={\bf
P}(({\wedge}^2W)^{\perp})\subset \check{Pf}$.
Then $\Pn 8_{W}$
intersects the Grassmannian $\Gr V2$
along the 5-fold Schubert cycle
$Y_{W} := {\sigma}_{30}(W^{\perp})$
   of $2$-dimensional subspaces
of $V^{*}$ that intersects the rank $2$ subspace $W^{\perp}$ nontrivially.
Therefore the general plane in $\Pn 8$
does not intersect $Y_{W}$.

When $A$ has rank $4$, the kernel is a rank $2$ subspace $U_{A}\subset
V$.  So there is a natural kernel map:
$$pker:\check{Pf}\setminus\Gr V2 \rightarrow \Gr 2V\quad [A] \mapsto
[U_{A}].$$
This map can also be seen as the map
$${\wedge}^2 V^{*}\to {\wedge}^4 V^{*}\cong
{\wedge}^2 V\quad \alpha\mapsto\alpha\wedge \alpha$$
so it is quadratic in the coordinates.  Now, the hyperplane
section $H_{A}\cap \Gr 2V$ is singular precisely in ${\bf
P}({\wedge}^2 U_{A})$.  On the other hand ${\bf
P}({\wedge}^2 U_{A})\in \Gr 2{W'}$ for any $4$-dimensional subspace $W'\subset
V$ that contains $U_{A}$.  Clearly $U_{A}\subset W$ for
any $A\in {\wedge}^2 W^{\perp}$.  Therefore
$$Z'=\cap_{A\in {\Pi}^2} H_{A}\cap \Gr 2V$$
contains the image of
$$s:{\Pi}^2\to \Gr 2V\quad A\mapsto {\bf
P}({\wedge}^2 U_{A}).$$
Since ${\bf P}({\wedge}^2 U_{A})$ is a singular point in
$H_{A}\cap\Gr 2V$, the image
of $s$ is contained in the singular locus of $Z'$.

\begin{lemma}\label{empty}
$Sing(Z')$ is contained in a Veronese surface if and only if  ${\Pi}^2
\cap \Gr V2 = \emptyset$.
\end{lemma}

\proof Indeed, the above map $s$ is defined everywhere on ${\Pi}^2$ only if
${\Pi}^2 \cap \Gr V2 = \emptyset$, and in this case the image
is clearly a Veronese surface.  On the other hand if $A \in {\Pi}^2 \cap
\Gr V2$,
and $U_{A} \subset V$
is the kernel of $A$ (regarded as a skew-symmetric
map as above), then the hyperplane section $H_A \subset \Gr 2V$
defined by $A$ is singular along a 4-fold quadric
$Q_A = \Gr 2{U_{A}}\subset \Gr 2V$.  Thus
$$Z'=\cap_{A\in {\Pi}^2} H_{A}\cap \Gr 2V$$
is singular at least along a
codimension 2 linear section of a $4-$fold quadric $Q_A$, which is
clearly not contained in a Veronese surface.
   \qed

\begin{proposition}\label{Zplane}
The linear section $Z = \overline{H}_\omega$
of $\Gr 2V$ is defined by
a plane ${\Pi}^2$ of linear equations in $\check{Pf} \setminus
\Gr V2$.  The singular locus of
$Z$ is a Veronese surface.

\end{proposition}

\proof  We noticed above that the
singular locus of $Z$ is contained in a Veronese surface.  It follows
from \ref{pfplane} and \ref{empty} that the singular locus of $Z$ is
a Veronese surface and the orthogonal plane does not intersect
$\Gr 2V^{*}$.
\qed

The following result on the orbits of the $GL(V)$-action on $\check{\bf P}^{14}$  fits well with the remark
above.
\begin{proposition}(Sato and Kimura 
\cite{SK} p. 94)\label{Zorbit}
The action ${\wedge}^2$
of $GL(V)$ on $\check{\bf P}^{14}$ is transitive
on the set of the planes
${\Pi}^2 \subset \check{Pf} \setminus \Gr V2$.
\end{proposition}
\ref{Zplane}, \ref{Zorbit} immediately implies the following 
corollary and the second part of theorem \ref{theoremA}.
\begin{corollary}\label{Zsurj}
    Let $Z\subset \Gr 2V$ be a $5$-fold linear section.  Then the 
    following are equivalent:\par\noindent
    i) $Z$ contains a 
$4$-dimensional smooth quadric and sing$Z$ is a Veronese 
surface.\par\noindent ii) The orthogonal complement of the span $<Z>^{\bot}\subset 
\check{Pf} \setminus
\Gr V2$.\par\noindent
iii) $Z$ is the projection of a nodal hyperplane section of $\Sigma=\LG 
36\subset\Pn {13}$ 
from its node.\end{corollary}
Let $S=Z\cap {\bf P}^{8}$ be a general linear surface section of 
$Z$.  Then, clearly, $S$ is a $K3$ surface with a conic section $C$.
\begin{corollary}\label{Zsurface}
The Picard group of a general linear surface section $S$ of $Z$ has 
rank $2$ and is generated by the class of a hyperplane and the class 
of the conic section on $S$.
\end{corollary}
\proof  According to the refined versions of Mukai's linear
section theorem \cite {M}, \cite {M4},  the general $K3$ surface $S$ with Picard group 
generated by a very ample line bundle ${\cal O}_{S}(H)$ of degree 
$H^2=14$ 
and the line bundle ${\cal O}_{S}(C)$ of a rational curve $C$ with 
$C\cdot H=2$ is a linear section of $\Gr 26$. The conic section $C$ 
lies then in a unique $4$-dimensional quadric $\Gr 24$ inside $\Gr 26$ 
and $S$ is therefore a linear section of a subvariety $Z$.\qed

\subsection{Stable rank 2 vector bundles on linear
sections}\label{stablevb}

\paragraph{}\label{la}
Theorem \ref{thm:hsection} allows us to construct families of rank two vector
bundles on linear sections of $\Sigma$ as promised in the end of
section \ref{subsec:projection}.
Let $1<k<6$ and consider
a smooth linear section $X=X^{6-k}=\Pn {13-k}\cap \Sigma$  and its
$Sp(3)$-dual linear section
$\check{F}(X)$ of the quartic $\check F$.
Then to each point  $\omega \in \check{F}(X)\setminus\check{\Omega}(X)$ 
we may associate a rank $2$ vector bundle $E_{\omega, X}$ on $X$
with chern classes $c_{1}(E_{\omega, X})=H$ and $c_{2}(E_{\omega, X})=\sigma_{X}$, where 
$\sigma_{X}$ is the class of a codimension $k-1$ linear section of a 
Segre $3$-fold.  In case $X$ is a curve $c_{2}(E_{\omega, X})=0$, but the vector 
bundle is special having at least $6$ global sections.

Namely, let $\omega \in \check{F}(X)\setminus\check{\Omega}(X)$,
then by \ref{thm:hsection}
   the projection ${\pi}_{u(\omega)}$ from the pivot $u(\omega) =
   L(piv^*(\omega))$
     sends the hyperplane section $H_\omega = \Sigma \cap {\bf
     P}^{12}_\omega$
      to the linear section
       $\overline{H}_\omega = \Gr 26 \cap \overline{\bf P}^{11}_\omega$.
Now $\omega \in  \check{F}(X) \subset ({\bf P}^{13-k})^{\perp}$,
so
$H_\omega \supset X$.
By assumption $X$ is smooth,
so \ref{singP} implies that
$u(\omega) \not\in {\bf P}^{13-k} = <X>$.
Therefore the projection
${\pi}_{u(\omega)} : {\bf P}^{12}_\omega\rightarrow \overline{\bf
P}^{11}_\omega$ restricts to a linear isomorphism of $X\subset\Pn
{13-k}$.
If $E_{\omega}$ is the rank 2 vector bundle on $H'_{\omega}$ constructed
in \ref{taut},
then the restriction
$E_{\omega, X}=E_\omega|_X$ is a rank 2 vector bundle on $X$.
Via the linear isomorphism $X\subset \Gr 26$, and $E_{\omega, X}$ is
the pullback of the universal rank 2 quotient bundle.  Therefore
$h^0(E_{\omega, X}) \ge 6$;
and  $c_{1}(E_{\omega, X}) = H_{X}$.  By \ref{taut}, the general section 
of $E_{\omega}$ vanishes on the projection of a Segre $3$-fold inside 
$H_{\omega}$
through $u(\omega)$.  Since $X$ does not pass through $u(\omega)$, the 
restriction to $X$ is that of a codimension $k-1$ linear section of 
this Segre $3$-fold, so $c_{2}(E_{\omega, X})=\sigma_{X}$. 

The vector bundle $E_{\omega, X}$ is in fact stable with respect to 
$H_{X}$.  Denote by  ${\cal M}_{X}(2,H,\sigma_{X})$ the moduli space of stable vector 
bundles on $X$ with chern classes $c_{1}(E)=H,c_{2}(E)=\sigma_{X}$.   
It exists as a quasiprojective variety (cf. \cite{Se}, \cite{Ma1}, \cite{Ma2}). 

\begin{proposition}\label{modmap} Let $1<k<6$, and let $X=X^{6-k}=
\Sigma\cap \Pn {13-k}$ be a smooth linear section
of $\Sigma$ without non-trivial automorphisms, and  let $\check{F}(X)$
be its $Sp(3)$-dual linear section of the quartic $\check F$.
Then there is a natural map
$$e_{X}: \check{F}(X)\setminus\check{\Omega}(X)\to {\cal M}_{X}(2,H,\sigma_{X})\quad \omega
\mapsto [E_{\omega, X}]$$
where $[E_{\omega, X}]$ is the isomorphism class of $E_{\omega, X}$, 
and $\sigma_{X}=0$ if $X$ is a curve.  Furthermore, the map is injective
on the set where $H^0(E_{\omega, X})=6$.
\end{proposition}
\proof  First we prove stability.
Recall that a rank $2$ vector bundle $E$ is  {\it Takemoto-Mumford stable (resp. semistable)} with
respect to the
   polarization $H$, if for each
line subbundle $L$, the inequality
$$2H^{i}\cdot L< H^{i}\cdot c_{1}(E), \quad ({\rm resp.} \quad
2H^{i}\cdot L\leq H^{i}\cdot c_{1}(E))$$
holds, where $i={\rm dim}X-1$. 

It is clear from the definion of Takemoto-Mumford that it is enough
in our case to check
stability in the curves $X^1$; the instability in the other cases
would imply instability by restriction to a curve section $X^1 
\subset X^{6-k}$.

So let $C=X^1\subset \Sigma $ be
a smooth linear curve section, and let  $\omega\in
\check{F}(C)\setminus \check{\Omega}(C)$.  We may assume that
$C$ has no automorphisms and that $C$ has no $g^1_{5}$  (see
\ref{gon} below).  Consider the vector 
bundle $E_{\omega,
C}$.  By assumption it is the pullback of the universal rank $2$ 
quotient bundle on $\Gr 26$, so the associated map of the ${\bf 
P}^1$-bundle $${\bf P}(E_{\omega, 
C})\to {\bf P}^5$$
is a morphism.  Let $S_{C}\subset{\bf P}^5$ be the image ruled 
surface of this 
morphism.  The curve $C$ is contained in a $5$-dimenional linear section $Z\subset 
\Gr 26$ that contains a $4$-dimensional quadric.  The linear span of $C$ 
intersects $Z$ in a surface, which by \ref{Zsurface} is a $K3$-surface section $Y$ containing a unique conic 
section that does not meet $C$. 
Assume now that $E_{\omega, 
C}$ has a subbundle of degree $d\geq 6$.   Then $S_{C}$ has $d$ members 
of the ruling contained in a hyperplane $H$.  In particular $C\subset \Gr 
26$ meets the Grassmannian $\Gr 2H$  in $d$ points.  But $\Gr 2H$ has 
degree $5$, so $\Gr 2H$ must intersect the surface $Y$ in at least a 
curve.   Since $Y$ is an irreducible surface, the intersection must be 
a curve and spans at most a $\Pn 4$.  Since $d\geq 6$, the corresponding divisor on $C$ of degree 
$d$ spans at least a $\Pn 4$.  Therefore this curve has degree $5$.  But 
by \ref{Zsurface} the Picard group of $Y$ is generated class of $H$ 
and the class of the unique conic section, so in particular every curve on it has even degree, a 
contradiction. 
Therefore $E_{\omega, 
C}$ is stable.

For injectivity, notice that
for a given embedding of $X\in \Gr 26$ the linear span of $X$ cuts the Grassmannian in
a variety $Y$ of
dimension dim$X +1$.  Given two elements $\omega$ and $\omega'$ in
$\check{F}(X)\setminus\check{\Omega}(X)$, then the
vector bundles
$E_{\omega, X}$, and $E_{\omega', X}$  are isomorphic only if the
two linear sections $Y_{\omega}$ and $Y_{\omega'}$ are projectively
equivalent.  In fact the global sections of $E_{\omega, X}$
define the
map into $\Gr 26$ and the linear span of the image defines $Y_{\omega}$. Now
$Y_{\omega}$ is the projection from $u(\omega)$ of the linear
section $\tilde Y_{\omega}$
of $\Sigma$ defined by the span of $X$ and $u(\omega)$.  Therefore
$Y_{\omega}$ and $Y_{\omega'}$ are projectively equivalent if and
only if $\tilde Y_{\omega}$ and $\tilde Y_{\omega'}$ are equivalent.  
Now, Mukai proves in \cite{M} Theorem 0.2, that two smooth linear surface sections 
of $\Sigma$ are projectively equivalent if and only if they lie in the 
same orbit of the group action $\rho$.  It is straightforward to extend his argument 
to our nodal case.  Therefore the linear span of $X$ and $u(\omega)$ and the linear
span of $X$ and $u(\omega') $ are in the same orbit under the action
$\rho$.  As soon as $X$ has no nontrivial automorphisms, this cannot
happen.
\qed

\medskip
We next see what is the image of the map $e_{X}$, and start with the 
curve case.
The general stable  rank $2$ vector bundle with canonical determinant 
on $C$ has no sections.
The subset of $M_C(2;K)$ corresponding to vector bundles with a given
number of sections has the structure of a subvariety, which have been
studied by several authors
(cf. \cite{M5}, \cite{M3}, \cite{OPP}, \cite{BF}).
Following their notation we define the Brill-Noether locus
$M_C(2;K,k)$ to be the subvariety of
$M_C(2;K)$ corresponding to vector bundles with at least $k+2$ sections.
\medskip
Let $E$ be a rank $2$ vector bundle on a general linear curve 
section $C\subset \Sigma$ with $[E]\in M_C(2;K,4)$.  Assume that $E$ 
is generated by global sections.  Let 
$$\wedge^2 H^0(C,E)\to H^0(C,K)$$
be the natural map.  If this map is surjective, then $E$ clearly is in 
the image of the map $e_{C}$ if the induced image $C\subset \Gr 
2{H^0(C,K)}$ is contained in a subvariety isomorphic to $Z$.  By 
\ref{Zsurj} this is 
satisfied as soon as 
the orthogonal complement of the span of $C$ in ${\bf P}(\wedge^2 
H^0(C,K))$ does not meet $\Gr 
{H^0(C,K)}2$. 
Thus
\begin{lemma} The isomorphism class $[E]\in M_C(2;K,4)$ fails to be in the image of $e_{C}$ 
only if 
$$\wedge^2 H^0(C,E)\to H^0(C,K)$$
is not surjective, or 
the orthogonal complement of the span of $C$ meets $\Gr 
4{H^0(C,K)}$, or $E$ is not globally generated.
\end{lemma}

We will not show that the cases of the lemma do not occur, but simply 
note that they represent closed subvarieties of $M_C(2;K,4)$.
A more general result is:

\begin{proposition}\label{bncurve} (\cite{M3} Theorem 4.7 or \cite{BF} p. 260-261):
Let $C$ be a curve of genus $9$ with no $g^1_{5}$.
If $M_C(2;K,5) = \emptyset$, then $M_C(2;K,4)$ is smooth and of dimension
three
precisely at the points representing bundles $E$ for which the Petri map
${\mu}: Sym^2 \ H^0(C,E) \rightarrow H^0(C,Sym^2 \ E)$
is injective.\end{proposition}

The injectivity of the Petri map is shown by Bertram and Feinberg for
$g(C)\geq 2$ and any stable rank $2$ vector bundle with canonical
determinant and $h^0(C, E)\leq 5$ in \cite{BF}.  The same line of
argument yields
\begin{lemma}\label{Petri}

Let $C$ be a curve of genus $9$ as above, then the Petri map

$$\mu :Sym^2 \ H^0(C,E) \rightarrow H^0(C,Sym^2 \ E)$$
\noindent
is injective for any stable bundle $E \in M_C(2;K,4)$.


\end{lemma}




\proof  Let $S$ be the scroll in ${\bf P}(H^0(C,E)^*)$ defined by mapping ${\bf 
P}(E)$ by the global sections of $E$. 
Then an element in the kernel of the Petri map 
$$\mu :Sym^2 \ H^0(C,E) \rightarrow H^0(C,Sym^2 \ E)$$
 defines a quadric hypersurface
containing $S$ (see \cite{BF} p. 267).
So the Petri map  is injective if and only if 
$S$ is contained in any quadric $Q$.  Since det$E$ is the canonical
line bundle, the degree of $S$ is 16. 
Since $C$ is a section of $\Sigma$ it has no $g^1_{5}, g^2_{7}$,
$g^3_{9}$ or $g^4_{11}$.  Furthermore, since $E$ is stable, it has no 
sub linebundle of degree greater than $7$, or equivalently, no section of $E$
vanishes in a divisor of degree $8$.  In particular no
$\Pn 4$ contains $8$ lines on $S$, and no plane intersects $S$ in a
section.  For the latter, it is clear that a plane section of degree
at most $7$ corresponds to a linear series on $C$ of dimension $2$
and degree at most $7$, while for a plane section of degree at least
$8$, the residual net of hyperplanes define a linear series of
dimension $2$ and degree at most $8$.  Equality corresponds to
a semistable bundle $E$.

Following \cite{BF} \S 4, we regard separately
the cases $1\leq rank \ Q \leq 6$:

\medskip

If $rank \ Q = 6$, then  
the family $F(Q) \subset \Gr 26$ of lines
on the smooth quadric $Q$
is isomorphic to the $5$-fold

$$F(Q) = {\bf P}(T_{\Pn 3}(-1))\subset \Pn 3 \times \Pnd 3, $$

\medskip

\noindent
which is the incidence variety
between points and planes in $\Pn 3$.
Denote the pullbacks of the hyperplane divisors in $\Pn 3$ and $\Pnd
3$ to $C$ by $h$ and $h'$.  Then $h$ and $h'$ are complimentary divisors
in a canonical divisor, i.e. $h+h'\equiv K_{C}$.

Assume that one of the two projections maps $C$ to a line.  This line
will then intersect all lines in $\Pn 3$ that are parameterized by
$S$.  But this means that $S$ is degenerate, contained in a $\Pn 4$,
a contradiction.

Thus the linear series defined by $h$ and $h'$ both have dimension at
least 2.  Since $C$ has no $g^2_{7}$, the degree of both $h$ and $h'$
is at least $8$.  Since they are complimentary in a canonical
divisor, this happens only if they both define $g^2_{8}$'s.
This corresponds to a semistable vector bundle $E$.
\medskip

If $rank \ Q = 5$,
then $Q$ is a cone
with vertex a point, and the planes in $Q$ all pass through the
vertex and are parameterized by $\Pn 3$.  Therefore $F(Q)\subset\Gr
26$ is a $\Pn
2$-bundle over a $\Pn 3$.  Since $S$ is no cone, only finitely many
lines of $S$, say $d$ lines, pass through the vertex of $Q$.  Let $P$
be the $\Pn 4$
of lines in $\Gr 26$ passing through the vertex, and let $p: C\to \Gr
25$, be the projection from $P$.  Then $p$ corresponds to the
projection of $S$ from the vertex of $Q$, and maps $C$ into the
double Veronese embedding of $\Pn 3$ in $\Gr 25$.  Thus the canonical
linear series has a decomposition as a sum $K_{C}=D+2L$, where
$D=C\cap P$ is a divisor of even degree $d$.  Since $D$ spans at most
a $\Pn 4$ in the canonical embedding of $C$, the degree $d\leq 6$.
If $d\leq 2$, then
$p(C)$ spans at least a $\Pn 6$, and $L$ is of degree at most $8$ and dimension
$3$ contrary to the assumption on $C$.  If $d= 4$ or $d=6$, then $L$ is
a $g^2_{6}$ resp. a $g^1_{5}$, again a contradiction.
\medskip

If $rank \ Q = 4$, then $Q$ has two pencils of $\Pn 3$'s.  The  restriction of these
pencils to $S$ define pencils of curves $|D|$ and $|D'|$ on $S$ such
that $D+D'$ is a hyperplane section.  We may assume that $D$ is a
section of $S$, while $D'$ is the pullback of a divisor on $C$.  Thus
deg$D'\geq 6$ and deg$D\leq 10$.  Since $C$ has no $g^1_{5}$ only equality is
possible.  In this case the decomposition $D+D'$ of a hyperplane section of $S$,
correspond to an exact sequence
$$0\to {\cal O}_{C}(D_{C}')\to E\to {\cal O}_{C}(K_{C}-D_{C}')\to 0.$$
The assumption on $C$ in fact implies that this sequence is exact on
global sections, i.e. the connecting homomorphism
${\delta}_{[E]}: H^0({\cal O}_{C}(K_{C}-D_{C}')) \rightarrow H^1({\cal
O}_C(D_{C}'))$
is zero.
Therefore $E = {\cal O}_{C}(D_{C}') \oplus {\cal O}_{C}(K_{C}-D_{C}')$;
and since $deg(K_{C}-D_{C}') = 10>8$ the bundle $E$ is not even
semistable.

\medskip

If $rank \ Q = 3$, then $Q$ is a cone with vertex a plane ${\bf P}^2$
over a smooth plane conic $q$.  Thus the hyperplane divisor on $S$
decomposes $H=2D+D_{0}$, where $D_{0}=S\cap \Pn 2$ is the intersection
of $S$ with the vertex of $Q$.   Note that $D_{0}$ has to be a curve, 
otherwise $H=2D$ contradicting the fact that $H$ is a section of $S$ 
over $C$.  Thus $D_{0}$ must
be a section of $S$, but we saw above that no section of $S$ lie in a plane 
so this case is impossible.
\medskip

If $rank \ Q = 2$ or $1$, then the scroll $S$ span at most a ${\bf P}^4$ contrary to
the assumption.
\qed

One part of Mukai's famous linear section theorem says:

\begin{theorem}\label{gon}  (Mukai  \cite{M}). {\sl Any smooth curve
$C$ of genus $9$ with no $g^1_5$
is isomorphic to a linear section $C=X^1$ of $\Sigma=\LG 36$. }
\end{theorem}

On the other hand 
\begin{lemma}\label{5gon} No smooth linear curve section $C$ in a 
$\Sigma$ has a $g^1_{5}$.
    \end{lemma}
    
    \proof Consider the curve $C$ as a subvariety of $\Gr 36$, and 
    let $D$ be a member of a $g^1_{5}$ on $C$.  Then $D$ spans a $\Pn 
    3$ and must correspond to five $3$-spaces in contained in a $\Gr 
    35$.  The intersection of $\Sigma$ with any $\Gr 35$ is however 
    always a $\Pn 4$ section of a Grassmannian quadric.  Therefore 
    the intersecton with the span of $D$ must be a quadric surface, a 
    contradiction with the fact that $C$ is a linear section.\qed

Therefore, if we combine \ref{Petri} and \ref{bncurve} and the injectivity of the map
$e_{C}$, we recover Mukai's result on the Brill-Noether locus:
\begin{theorem}\label{curvecomp}(\cite{M4} p.17)  For a smooth linear curve 
section $C$ of
$\Sigma$ the quartic $3$-fold $\check{F}(C)\setminus \check{\Omega}(C)$ is a 
connected component of
the Brill-Noether locus $M_C(2;K,4)$. The $21$ double
points $\check{\Omega}(C)\subset \check{F}(C)$ in the 
boundary correspond to
semistable vector bundles that are not stable.
\end{theorem}

\proof It only
remains to check the semistable
boundary.
The {\it semistable boundary} ${\delta}_{ss}M_C(2;K)$
of $M_C(2;K)$
is the image of $Pic^{g-1}(C)$ under the map

$$j: Pic^{g-1}(C) \rightarrow M_C(2;K),
              \ j: L \mapsto L\oplus K \otimes L^{-1}$$

\noindent
-- see \cite{OPP} \S 1.
The {\it semistable boundary} of the locus
$M_C(2;K,4) \subset M_C(2;K)$ is the intersection
${\delta}_{ss}M_C(2;K,4) = {\delta}_{ss}M_C(2;K) \cap M_C(2;K,4)$.
Therefore

$${\delta}_{ss}M_C(2;K,4)
=
\{ L \in Pic^{8}(C): L \oplus K \otimes L^{-1} \in M_C(2;K,4)$$
$$=\{ L \in Pic^{8}(C): h^0(L \oplus K \otimes L^{-1}) \ge 6 \}.$$
Since $C\subset \LG 36$ has no  $g^1_{5}$ it has no $g^3_{8}$, so any 
line bundle
$L$ and likewise $K\otimes 
L^{-1}$, such that $h^0(L \oplus K \otimes 
L^{-1}) \ge 6$ must be a $g^2_{8}$.
Let $W^r_d(C) \subset Pic^d(C)$
be the Brill-Noether locus of all the invertible sheaves
$L$ of degree $d$ on $C$ such that $h^0(C,L) \ge r$.
Since $C$ is general of genus $g = 9$ then the fundamental class
of $W^r_d(C)$ in $Pic^d(C) \cong J(C)$ is

$$[W^r_d]
=
\frac{r!.(r-1)!...0!}{(g+2r-d)!...(g+r-d)!}{\Theta}^{(r+1)(g+r-d)},$$
\noindent
where $(J(C),{\Theta})$ is the principally
polarized jacobian of $C$
(see \cite{GH}, Ch. 2 \S 7 -- Special linear systems IV).
In particular $dim \ W^2_8(C) = 0$;
and since $deg({\Theta}^9/9!) = 1$ 

$$deg \ W^2_8(C) = \frac{2!.1.!0!}{5!.4!.3!}9! = 42.$$

Therefore on the general curve $C$ of genus $9$ there are
exactly $42$ line bundles $L$ such that $deg(L) = 8$ and
$h^0(C,L) = 3$.
Moreover, since $K\otimes 
L^{-1}$ also has degree $8$ and $3$ sections the map

$$\bar{}: W^2_8(C) \rightarrow W^2_8(C),\quad
\bar{}: L \mapsto \bar{L} = K \otimes L^{-1}$$

\noindent
is an involution of $W^2_8(C)$.
The fixed points, if such exist,
of the involution \ $\bar{}$ \ are
these $L$ such that $L^{\otimes 2} = K_C$
(i.e. $L$ is a theta-characteristic of $C$)
for which $h^0(C,L) = 3$.
But since $C$ is general $h^0(C,L) \le 1$ for
any theta-characteristic of $C$, i.e. \ $\bar{}$ \
has no fixed points.

Therefore on the general curve $C$ of genus $9$ there exist
exactly $21$ (non-ordered) pairs $(L_i,\bar{L}_i), 1 \le i \le 21$
of line bundles such that \
$deg \ L_i = deg \ \bar{L}_i = 8$,
$h^0(C,L_i) = h^0(C,\bar{L}_i) = 3$
and $L_i \otimes \bar{L}_i = K_C$.

Therefore the semistable
boundary ${\delta}_{ss}M_C(2;K,4)$ of $M_C(2;K,4)$
is a finite set of $21$ points representing the $21$ rank $2$
vector bundles $E_i = L_i \oplus \bar{L}_i, 1 \le i \le 21$.

                                                                         %
In our setting, when $C$ is a general linear section of $\Sigma$, the
$Sp(3)$-dual $\check{F}(C)$ is a quartic $3$-fold that intersects
$\check\Omega$ in $21$ nodes. Since this number fits with the number 
of semistable vector bundles just computed we try to extend the map
$e_{C}$ above to $\check{\Omega}(C)$.  This is possible:  Let $\omega
\in \check\Omega$ and let $u$ be a dual pivot of $\omega$ and consider
the blowup of $\Sigma$ in $u$.  Let $Q'_{\omega}$ be the strict
transform of the exceptional divisor on $H'_{\omega}$, as above.
Consider the exact sequence
$$0\to E'_{x}\to \wedge^2 U^{*}_{\omega}(-Q'_{\omega})\to {\cal
O}_{H'_{\omega}}(L')\to N_{x}\to 0,$$
where the cokernel sheaf $N_{x}={\cal
O}_{C(x)}(L')$ is the restriction of the
line bundle ${\cal
O}_{\Sigma_{\omega}}(L')$ to the zero-scheme of
$x$.  By $\ref{nonsurjective}$ the kernel sheaf $E'_{x}$ has $6$
sections as soon as $\Pn 4_{v}\subset \Pn {12}_{\omega}$, where
$v=L^{-1}(x)$.  So $v$ belongs to one of the planes in the involutive pair
$P_{1}$ and $P_{2}$ belonging to $\omega$.  Furthermore the zeros of
$x$ have codimension
$2$ so it does not intersect a general curve section $C\subset
H_{\omega}$.  Therefore the restriction of $E'_{x}$ to $C$ becomes a
rank $2$ vector bundle $E_{\omega}$ with canonical determinant and
$6$ sections.  On the other hand for any line $l$ in one of the planes
$P_{i}$, the subvariety $\Sigma_{l}$ of Lagrangian planes that meet
$l$, is a Weil-divisor on $H_{\omega}$.  Thus there are two nets of
Weil divisors on $H_{\omega}$. Restricted to the curve $C$ these
divisors become Cartier divisors, and the vector bundle $E_{\omega}$
splits as the sum of the corresponding line bundles. This yields the 
desired semistable vector bundle corresponding to the point 
$\omega\in \check\Omega (C)$. \qed

\bigskip

In the next section we study the image of the map $e_{X}$ of
\ref{modmap} in
the cases where $X$ is a surface, a threefold or a fourfold.  Thus we
recover and generalize the results by Mukai that initiated this
investigation.

\subsection{The moduli space $M_S(2;h,4)$ for a K3-surface $S$
of genus $9$}

\paragraph{}\label{ma}
Let $(S,h)$ be a polarized $K3$ surface of genus $g = 2n+1$,
let $s$ be an integer such that $s \le n$.
By \S 10 of \cite{M84} or \S 3 of \cite{M88},
the moduli space of stable rank 2 vector bundles $E$ on $S$:

$$M_S(2,h,s) = \{ E | E \mbox{ is stable, }
c_1(E) = h \mbox { and } \chi(S,E) = s+2 \}/(iso).
$$

\noindent
is a nonsingular symplectic variety of dimension $2(g-2s)$.
In particular if $s = n = (g-1)/2$ then
$\hat{S} := M_S(2,h,s)$
is a $K3$ surface.  Irreducibility has been given several proofs, see  
\cite{Go} for a recent one and including references to other proofs.

\medskip

Let $S = S_{16}$ be a general $K3$ surface of genus $9$ embedded
as a linear section of $\Sigma$
by a codimension $4$ subspace ${\bf P}^9 \subset {\bf P}^{13}$,
and let $H$ be the hyperplane class of $S \subset {\bf P}^9$.
The moduli space $M_S(S,H,\sigma_{S})$, then coincides with 
$M_S(S,H,4)$ since $\chi(S,E) = 4+2 = c_{2}(E)= {\rm deg}\sigma_{S}$.  
Therefore the above combine with \ref{modmap} to yield

\begin{theorem}\label{bnosurface}
For the general linear surface section $S=X\subset \Sigma$ the $K3$
surface $\hat{S} = M_S(2,H,4)$
is isomorphic to
the $Sp(3)$-dual quartic surface $\check{F}(S)$.
\end{theorem}
\proof
The map $e_{X}$ is injective, so $\check{F}(S)$ is a subvariety of 
$M_S(2,H,4)$.  On the other hand $M_S(2,H,4)$ is a $K3$ surface so 
they must coincide.
\qed

If we compare this with \ref{modmap} we see that in fact the map 
$e_{S}$ for linear surface sections $S\subset\Sigma$ is surjective.

\begin{proposition}\label{surjective}
    Let $X$ be a general smooth threefold or fourfold linear section
of $\Sigma$ and let $E$ be a stable
rank two vector bundle on $X$
with $h^0(X,E)=6$ and det$E={\cal O}_{X}(H)$.  Assume that natural
map $\wedge^2 H^0(X,E)\to
H^0({\cal O}_{X}(H))$ is surjective.  Then $E$ is in the image of $e_{X}$.
\end{proposition}

\proof  Assume that $E$ is a rank $2$ vector bundle on $X$ that
satisfies the conditions of the proposition.  Then the surjection
$\wedge^2 H^0(X,E)\to
H^0({\cal O}_{X}(H))$ defines a embedding $X\subset \Gr
2{H^0(X,E)^{*}}$.  Then clearly $E$ is in the image of $e_{X}$ if and
only if there is a $\Pn {11}$ such that $X\subset Z_{X}=\Pn {11}\cap \Gr
2{H^0(X,E)^{*}} $ for some $Z_{X}\cong Z$ of \ref{Zplane}.

Assume that dim$X\geq 3$. Since $e_{S}$ is surjective for any general 
surface section each surface section $S$ intersects the Grassmannian in a
threefold that is a linear section of a variety $Z_{S}\cong Z$.
Therefore the $4$-dimensional quadric $Q\subset Z_{S}$ intersects the linear
span of $S$ in a quadric surface.
If two surface sections $S$ and $S'$ of $X$ give rise to
subvarieties $Z_{S}$ and $Z_{S'}$ with distinct quadrics $Q$ and
$Q'$, then these
quadrics are Grassmannians $\Gr
2W$ and  $\Gr
2W'$ for $4$-dimensional subspaces $W$ and $W'$ of $H^0(X,E)^{*}$. 
So $Q$ and $Q'$
have a plane or a point in common.  But $S$ and $S'$ are linear 
sections of $X$ and the corresponding quadric surfaces may be chosen to
be smooth with exactly a conic section in common.  This is a 
contradiction.   Therefore the
subvarieties $Z$ for distinct surface sections of $X$ have the same
$4$-dimensional quadric $Q$.

The linear span of $X$ intersects the quadric $Q$ in a linear section of dimension
dim$X-1$ of the quadric.  Therefore the linear span of $X$ and $Q$ has dimension
$11$ and cuts $\Gr 26$ in a subvariety projectively  equivalent to
$Z$.  \qed

\paragraph{}\label{na}
By \cite{M2} any smooth Fano 3-fold $X = X_{16}$
of degree $16$, of $rank \ Pic = 1$ and of index $1$
(or -- simply -- a {\it prime} Fano 3-fold of degree $16$)
is a linear $3$-fold section of
$\Sigma $, and the hyperplane class $H$ of
$X_{16} = \Sigma \cap {\bf P}^{10}$ is the ample generator
of $Pic \ X$ over ${\bf Z}$.

If $X = X_{16}$ be general, then the $Sp(3)$-dual linear section $\check{F}(X)$
is a smooth plane quartic curve which does not intersect
$\check{\Omega}$.

\begin{proposition}\label{bnothreefold}
Let $X = X_{16} \subset {\bf P}^{10}$ be a general
prime Fano threefold of degree $16$.
Then the $Sp(3)$-dual to $X$ plane quartic curve $\check{F}(X)$
is isomorphic to an irreducible component of the moduli
space $M_{X}(2;H,\sigma_{X})$
of stable rank $2$ vector bundles on $X$ with $c_1 = [h]$
and $c_2 = \sigma_{X}$, where $[h]$ is the class of the hyperplane
section and $\sigma_{X}$ is the class of a sextic elliptic curve on $X$.
\end{proposition}

\proof  
Now the condition in 
\ref{surjective} is certainly an open one, so the complement of the image 
of the map $e_{X}$ is closed.  Since $e_{X}$ is injective by \ref{modmap} and its image is 
closed the theorem follows.\qed 
\medskip

We end with an easy corollary in the fourfold case:
\begin{corollary}\label{bnofourfold}
   For a general linear fourfold section $X\subset\Sigma$
the $Sp(3)$-dual $\check{F}(X)$, consists of four points.  Whenever
$X$ has no automorphisms, these four points define precisely the four isomorphism
classes of stable rank 2 vector bundles $E$ with $c_{1}(E)=H$ and 
$c_{2}(E)=\sigma_{X}$ where $\sigma_{X}$ is the class of a Del Pezzo 
surface of degree $6$ on $X$ such that the natural map $\wedge^2 H^0(X,E)\to
H^0({\cal O}_{X}(H))$ is surjective.
\end{corollary}

\smallskip

{\small
${}$\hspace{1.2cm}{\bf Atanas Iliev}
\hspace{5.0cm}{\bf Kristian Ranestad}

${}$\hspace{1.2cm}Institute of Mathematics
\hspace{3.2cm}Matematisk Institutt, UiO

${}$\hspace{1.2cm}Bulgarian Academy of Sciences
\hspace{2.1cm}P.B. 1053  Blindern

${}$\hspace{1.2cm}Acad. G. Bonchev Str., 8
\hspace{3.1cm}N-0316 Oslo, Norway

${}$\hspace{1.2cm}1113 Sofia, Bulgaria
\hspace{4.0cm}e-mail: ranestad@math.uio.no

${}$\hspace{1.2cm}e-mail: ailiev@math.bas.bg
\hspace{2.7cm}
}

\end{document}